\documentclass[a4paper,11pt]{amsart}

\usepackage{amssymb,amsmath}%,eufrak}
\author[Florent Benaych-Georges]{Florent Benaych-Georges}\address{LPMA, UPMC Univ Paris 6, Case courier 188, 4, Place Jussieu, 75252 Paris Cedex 05, France. } \email{florent.benaych@gmail.com}

\title[Rectangular random matrices, related convolution]{Rectangular random matrices, related convolution}
\date{\today}
\newcommand{\bck}{\backslash}
\newcommand{\uo}{\underline{\textrm{Case $1^\circ$)}} }
\newcommand{\deuzio}{\underline{\textrm{Case $2^\circ$)}} }
\newcommand{\udeuzio}{\underline{\textrm{Cases $1^\circ$) and  $2^\circ$)}} }
\newcommand{\uov}{$1^\circ$), }
\newcommand{\uop}{$1^\circ$). }

\newcommand{\NC}{\operatorname{NC}}
\newcommand{\NCp}{\operatorname{NC'}}
\newcommand{\NCd}{\operatorname{NC_d}}

\newcommand{\ED}{\operatorname{E}}

\newcommand{\ds}{\displaystyle}

\newcommand{\Pro}{\operatorname{P}}

\newcommand{\tr}{\operatorname{tr}}
\newcommand{\rg}{\operatorname{rg}}

\newcommand{\mf}{\mathfrak}

\newcommand{\Tr}{\operatorname{Tr}}
\newcommand{\ninf}{\underset{n\to\infty}{\longrightarrow}}
\newcommand{\ssi}{if and only if }

\newcommand{\teo}{theorem }
\newcommand{\teov}{theorem, }

\newcommand{\fid}{$\sst\boxplus$-infinitely divisible distribution }

\newcommand{\E}{\mathbb{E}}

\newcommand{\R}{\mathbb{R}}
\newcommand{\C}{\mathbb{C}}
\newcommand{\n}{\mathbb{N}}
\newcommand{\z}{\mathbb{Z}}

\newcommand{\ud}{\mathrm{d}}

\newcommand{\trv}{Voiculescu transform }

\newcommand{\pro}{probability }

\newcommand{\f}{\frac}
\newcommand{\ff}{\frac{1}}
\newcommand{\lf}{\left}
\newcommand{\ri}{\right}
\newcommand{\spl}{spectral law }
\newcommand{\st}{such that }
\newcommand{\la}{\lambda}
\newcommand{\La}{\Lambda}
\newcommand{\tii}{\scriptstyle\times\ds\!}

\newcommand{\vfi}{\varphi}
\newcommand{\ste}{\, ;\, }
\newcommand{\mc}{\mathcal }

\newcommand{\eps}{\varepsilon}
\newcommand{\arc}{{\scriptscriptstyle\boxplus_{\la}}}
\newcommand{\arco}{{\scriptscriptstyle\boxplus_{0}}}
\newcommand{\gab}{\Delta_{\alpha,\beta}}
\newcommand{\sst}{\scriptstyle}
\newcommand{\kk}{\mathfrak{K}}
\newcommand{\bxp}{{\scriptscriptstyle\boxplus}}
\newcommand{\bxt}{{\scriptscriptstyle\boxtimes}}

\newcommand{\A}{\mc{A}}
\newcommand{\D}{\mc{D}}

\newtheorem{Th}{Theorem}[section]
\newtheorem{propo}[Th]{Proposition} 
 
\newtheorem{lem}[Th]{Lemma}
\newtheorem{ex}[Th]{Example}
\newtheorem{rmq}[Th]{Remark}
\newtheorem{cor}[Th]{Corollary}
\newtheorem{Def}[Th]{Definition}
\newtheorem{propdef}[Th]{Proposition-Definition}
\newenvironment{pr}{\noindent {\bf Proof. }}{\ \ \ $\square$}
\newenvironment{prth}{\noindent {\bf Proof of the theorem. }}{\ \ \ $\square$}

\long\def\symbolfootnote[#1]#2{\begingroup
\def\thefootnote{\fnsymbol{footnote}}\footnote[#1]{#2}\endgroup}

\begin{document}
\maketitle
\symbolfootnote[0]{{\it MSC 2000 subject classifications.}  15A52, %Random matrices 
46L54} %Free proabability and free factors

\symbolfootnote[0]{{\it Key words.} random matrices, free probability, free convolution} 

\begin{abstract}We characterize asymptotic collective behavior of rectangular random matrices, the sizes of which tend to infinity at different rates. It appears that one can compute the limits of all non commutative moments (thus all spectral properties) of the random matrices we consider because, when embedded in a space of larger square matrices, independent rectangular random matrices are asymptotically free with amalgamation over a subalgebra. Therefore, we can define a ``rectangular free convolution'', which allows to deduce the singular values of the sum of two large independent rectangular random matrices from the individual singular values. This convolution is linearized by cumulants and by an analytic integral transform, that we called the ``rectangular $R$-transform''. 
\end{abstract}

%\tableofcontents

\section*{Introduction} 
The first problem we are going to deal with in this paper is the modeling of  asymptotic collective behavior of independent rectangular random matrices. In order to explain to the reader the way we will treat this problem, let us recall him the work already done by Wigner, Pastur, Marchenko, Girko, Bai, Voiculescu,... for square random matrices. First of all, in the 50's Wigner considered self-adjoint random matrices with gaussian entries (GUE) and proved that the spectral law (i.e. the uniform distribution on the set of eigenvalues) of these random matrices %, which is a random probability law, 
converges %weakly in probability (even almost surely) 
to the so-called {\it semicircle law.}  This result was improved, and other results giving the asymptotic spectral law of random matrices were proved (see, among many other sources, \cite{pasturlejay}). In the same time, people studied the local structure of the spectrum of random matrices (see, e.g. \cite{mehta}), but this is not the kind of problem we are going to study here. A new point of view was adopted in the early 90's by Voiculescu, who proposed a way to compute the asymptotic normalized trace (when the dimension $n$ of the matrices goes to infinity) of products $$M(s_1,n)^{\eps_1}\cdots M(s_k,n)^{\eps_k}$$ of random matrices taken among a family $M(1,n),M(2,n),\ldots$ of independent random $n\tii n$ matrices and their adjoints, the only hypothesis being a bound on the norms of the matrices, the fact that the matrices   $M(1,n),M(2,n),\ldots$ have all a limit singular law (the {\it singular law} of a rectangular matrix is the uniform distribution on the set of its singular values, i.e. of the eigenvalues of the absolute value of the matrix), and  a property of invariance  of the distributions under an action of the unitary group. These results can be deduced from the fundamental article \cite{voic91}, but are presented under this form in \cite{hiai}.  The advantage of being able to compute the limit of such normalized traces is that it gives us the asymptotic normalized trace of any noncommutative polynomial in our independent random matrices. Hence, since the normalized trace of the $k$-th power of a matrix is the $k$-th moment of its spectral law,   we are able to give the asymptotic singular law of any polynomial of our random matrices. This is why this work is said to {\it model the asymptotic collective behavior of independent square random matrices.} 
For example, it can be proved (combine results of \cite{hiai} and \cite{haag2}) that the  asymptotic singular law of the sum of two independent random matrices whose distributions are invariant under left and right actions of the unitary group and whose  singular laws converge weakly to \pro measures $\mu_1$, $\mu_2$, only depends on $\mu_1$ and $\mu_2$, and can be expressed easily from $\mu_1$ and $\mu_2$: it is the \pro measure on $[0,\infty)$, the symmetrization of which is the free convolution of the symmetrizations of $\mu_1$, $\mu_2$. We often present the similar result for hermitian matrices, for which we work with spectral law in the place of singular  law, and for which no  symmetrization is necessary, but here, we shall work with rectangular matrices, which cannot be hermitian, so the square analogue of our work will be found in non hermitian matrices. In this text, as in the work of Voiculescu presented above, we will propose a way to compute the asymptotic normalized trace of  products of random matrices taken among a family of independent rectangular random matrices, whose sizes tend to infinity, but with different rates. The notion involved, similarly to freeness in Voiculescu's modeling of asymptotics of square random matrices, is freeness with amalgamation over a finite dimensional subalgebra. This notion   arises from operator-valued free \pro theory, but we chose to use the point of view of operator-valued free probability only when necessary, %because the problems we are going to deal with in this paper are of probabilistic essence, and 
because this point of view is not  satisfying in all cases: if one uses this point of view, then he has to consider the case where the ratio of sizes of our random matrices tends to zero separately. Moreover,  the point of view of operator-valued free \pro in the asymptotics of rectangular random matrices is developed in another paper, where we analyze the related free entropy and Fischer information (\cite{fbg.free.amalg}).                      

This modeling of asymptotics of rectangular random matrices  will allow us to define, for $\la\in [0,1]$, a binary operation $\arc$ on the set of symmetric \pro measures, called free convolution with ratio $\la$, and denoted by $\arc$.  For $\mu_1,\mu_2$ symmetric \pro measures, $\mu_1\arc \mu_2$ is defined to be the limit of the singular law  of a sum of two  independent rectangular random matrices, whose dimensions tend to infinity in a ratio $\la$, one of them being bi-unitarily invariant, and whose singular laws tend to $\mu_1,\mu_2$.  
In the second part of this paper, after having analyzed related cumulants, we construct an analytic integral transform
 which linearizes  $\arc$ (like the $R$-transform does for free convolution) and give examples. This part of this paper is the base of other researches. Firstly, in \cite{fbg05.inf.div}, we study the related infinite divisibility: it is proved that the set of $\arc$-infinitely divisible distributions is in a deep  correspondence with the  set of symmetric classical infinitely divisible distributions. Secondly, in \cite{BBG07}, we study  some questions related to the support and the regularity of $\mu_1\arc \mu_2$.

{\bf Acknowledgments.} We would like to thank Philippe Biane, Dan Voiculescu, and Piotr \'Sniady for useful discussions, as well as Thierry Cabanal-Duvillard, who organized the workshop ``Journ\'ee Probabilit\'es Libres'' at MAP5 in June 2004, where the author had the opportunity to have some of these discussions. Also, we would like to thank C\'ecile Martineau for her contribution to the english version of this paper.

\section{Asymptotic behavior of rectangular random matrices}
\subsection{Notes on reduction of rectangular random matrices}\label{14.INCENDIE} In this subsection, we recall some simple facts about polar decomposition of rectangular complex matrices, which can be found in \cite{H&J}. Consider a $p\tii q$ matrix $M$.  
When $ p\leq q$, we will denote by $|M|$ the only $p\tii p$ positive hermitian matrix $H$ \st we can write $M=HT$, with $T$ a $p\tii q$   matrix \st $TT^*=I_p$. In this case, the matrix $|M|$ is the square root of $MM^*$. When $p>q$, we will denote by $|M|$ the only $q\tii q$ positive hermitian matrix $H$ \st we can write $M=TH$, with $T$ a $p\tii q$   matrix \st $T^*T=I_q$.  In this case, the matrix $|M|$ is the square root of $M^*M$. 
In both cases, the spectrum of $|M|$ is the only (up to a permutation) family $(h_1,h_2\ldots)$ of nonnegative real numbers  \st one can write $$M=U[\delta_i^jh_i]_{\substack{ 1\leq i \leq p\\ 1\leq j\leq q}}V,$$
with $U$,  $V$ respectively $p\tii p$, $q\tii q$ unitary matrices. These numbers are called {\it singular values} of $M$. The uniform distribution on $h_1,h_2$... will be called the {\it singular law} of $M$ (whereas the uniform distribution on the spectrum of an hermitian matrix is called its {\it spectral law}). Note that in both cases, for any $\alpha >0$, $$\ff{\min(p,q)}\Tr |M|^\alpha =\ff{\min(p,q)}\Tr (MM^*)^{\alpha/\!2}=\ff{\min(p,q)}\Tr (M^*M)^{\alpha/\!2}$$ is the $\alpha$-th moment of the singular law of $M$.

A random matrix is said to be {\it bi-unitarily invariant} if its  distribution is invariant under the left and right actions of the unitary group.

\subsection{Definitions}\label{23.1.05.132} For definitions of  {\it singular values, singular law, spectral law, (bi-)unitarily invariant and uniform random unitary matrices}, we refer to subsection \ref{14.INCENDIE}. For all positive integer $d$, we denote the set $\{1,\ldots, d\}$ by $[d]$. We denote the normalized trace of a square matrix $X$ by $\tr X$ (whatever the size of $X$ is).

Consider  a positive integer $d$, and sequences $q_1(n),\ldots,q_d(n)$ of pairwise distinct positive integers which all tend to infinity as $n$ tends to infinity, and such that for all $k\in [d]$, \begin{equation}\label{22.3.05.1}\f{q_k(n)}{n}\ninf \rho_k\geq 0,\end{equation} at most one of the $\rho_k$'s being zero. We are going to work with rectangular random matrices of sizes $q_k(n)\tii q_l(n)$.  Free probability's modeling of the asymptotic behavior of square random matrices relies on a comparison between random matrices and elements of an algebra arising from operator algebra theory. To product analogous results for rectangular random matrices, we have embedded our matrices of sizes $q_k(n)\tii q_l(n)$ ($k,l\in [d]$) in an algebra.
 Let us assume that for all $n$,   $q_1(n)+\cdots +q_d(n)=n$ (since what we will prove would obviously also work when $n$ is replaced by a subsequence, it is not a real restriction). Now, $n\tii n$ matrices will be represented as $d^2$ block matrices, \st for all $k,l\in[d]$, the $(k,l)$-th block is a $q_k(n)\tii q_l(n)$ matrix. 

For all $k,l\in [d]$, for all $q_k(n)\tii q_l(n)$ matrix $M$, let us denote by $\widetilde{M}$ the ``$n\tii n$ extension of $M$'', that is the $n\tii n$ matrix   with $(i,j)$-th block $M$ if $(i,j)=(k,l)$, and zero in the other case. Note we have preservation of adjoints ($\widetilde{M^*}=\widetilde{M}^*$), and of products: $$\forall k,l,k',l'\in [d], \forall M\textrm{ of size $q_k(n)\tii q_l(n)$}, \forall N\textrm{ of size $q_{k'}(n)\tii q_{l'}(n)$}, \widetilde{M}\widetilde{N}=\begin{cases}\widetilde{MN}&\textrm{if $l=k'$,}\\ 0&\textrm{in the other case.}\end{cases}$$
%Moreover, the trace of $\widetilde{M}$ is the one of $M$ if  $M$ is square, and zero in the other case. Since we are overall interested in singular laws, whose moments are given by normalized  traces of square matrices, these embeddings shall not cause any loss of information. Note however that the normalized trace is not exactly preserved by these embeddings (the normalized trace of a matrix $M$ of size $q_k(n)\tii q_k(n)$ is $\f{n}{q_k(n)}$ times the normalized trace of $\widetilde{M}$), and that the well known relation $\Tr MN=\Tr NM$ is not satisfied anymore when we work with normalized traces of rectangular matrices. Indeed,  for all $k,l\in [d]$, for all $M$ of size $q_k(n)\tii q_l(n)$, for all $N$ of size $q_{l}(n)\tii q_{k}(n)$, we have  \begin{equation}\label{10.1.05.1} q_k(n)\tr MN=q_l(n) \tr NM.\end{equation}
 At last, let us define, for all $k\in [d]$, the projector $$p_k(n):=\widetilde{I_{q_k(n)}}.$$

 Let us formalize the structure we inherit. Consider a $*$-algebra  $\A$ endowed with a family $(p_1,\ldots,p_d)$ of non zero self-adjoint projectors (i.e. $\forall i, p_i^2=p_i$) which are pairwise orthogonal (i.e. $\forall i\neq j, p_ip_j=0$), and such that $p_1+\cdots +p_d=1$. Any element $x$ of $\mc{A}$ can then be represented $$x=\begin{bmatrix}x_{11}& \cdots & x_{1d}\\  \vdots & & \vdots \\ x_{d1}& \cdots & x_{dd}\end{bmatrix},$$ where $\forall i,j, x_{ij}=p_ixp_j$. This notation is compatible with the product and the involution. 
Assume each subalgebra $p_k\A p_k$ to be endowed with a tracial state $\vfi_k$ (i.e. $\vfi_k(p_k)=1$ and for all $x,y\in p_k\A p_k, \vfi_k(xy-yx)=0$) such that for all $ k,l\in [d], x\in p_k\A p_l,y\in p_l\A p_k, $ 
\begin{equation}\label{10.1.05.2}\rho_k\vfi_k(xy)=\rho_l\vfi_l(yx), \end{equation} where $(\rho_1,\ldots, \rho_d)$ is still the sequence of nonnegative real numbers defined by (\ref{22.3.05.1}). Recall that at most one of the $\rho_k$'s is zero. 
\begin{Def} Such a family $(\A,p_1,\ldots,p_d,\vfi_1,\ldots,\vfi_d)$ will be called a {\rm $(\rho_1,\ldots,\rho_d)$-rectangular probability space}. Elements  of the union of the $p_k\A p_l$'s ($k,l\in [d]$) will be called {\rm simple elements}.\end{Def}

\begin{ex}\label{camille.tango.triste}$\A$ is the $*$-algebra of $n\tii n$ complex matrices, $p_1,\ldots,p_d$ are $p_1(n),\ldots,p_d(n)$ previously defined, and each $\vfi_k$ is $\ff{q_k(n)}\Tr$. Then %by (\ref{10.1.05.1}), 
(\ref{10.1.05.2}) is satisfied when each $\rho_k$ is replaced by $q_k(n)/\! n$. \end{ex} 

%\begin{ex} If $\A$ is a $*$-algebra endowed with an orthogonal family of self-adjoint projectors $(p_1,\ldots, p_d)$ with sum $1$, if $\A$ is endowed with a tracial state $\vfi$ \st for all $k\in [d]$, $\rho_k:=\vfi(p_k)\neq 0$, then the family  $(\A,p_1,\ldots,p_d,\vfi_1,\ldots,\vfi_d)$ is a $(\rho_1,\ldots,\rho_d)$-rectangular probability space, where for all $k$, $\vfi_k=\ff{\rho_k}\vfi_{|p_k\A p_k}$. \end{ex} 

Consider a $(\rho_1,\ldots,\rho_d)$-rectangular probability space $(\A,p_1,\ldots,p_d,\vfi_1,\ldots,\vfi_d)$. Denote by $\D$ the linear span of the $p_k$'s, then $\D$ is an algebra %, which can be identified to the the set of $d\tii d$ complex diagonal matrices, and then to $\C^d$ by $$\sum_{k=1}^d\la_k p_k\simeq \diag (\la_1,\ldots,\la_d)\simeq (\la_1,\ldots,\la_d).$$ 
and the function $\ED: \A\to\C^d$, which maps $x\in \mc{A}$ to $\ED(x)=\sum_{k=1}^d\vfi_k(x_{kk}) p_k$ is a  {\it conditional expectation} from $\mc{A}$ to $\D$, i.e. $\ED(1)=1$ and  $\forall (d,a,d')\in \mc{D}\times \mc{A}\times \mc{D}, \ED(dad')=d\ED(a)d'.$

The following definition, first appeared in  \cite{voic95}, gives the right notion to describe asymptotics of independent random matrices.
\begin{Def}In an $*$-algebra $\mc{B}$ endowed with a $*$-subalgebra $\mc{C}$ and  a conditional  expectation $\ED_{\mc{C}}$ from $\mc{B}$ to $\mc{C}$, a family $(\mc{C}\subset \mc{B}_\alpha)_{\alpha\in A}$ of $*$-subalgebras is said to be {\rm free with amalgamation over $\mc{C}$} if for all $m\geq 1$, for all $\alpha_1\neq \cdots \neq \alpha_m\in A$, for all $x_1,\ldots,x_m$ elements of respectively $\mc{B}_{\alpha_1}$, \ldots,  $\mc{B}_{\alpha_m}$ one has  \begin{equation}\label{10.1.05.3.22.3.05}\ED_{\mc{C}}(x_1)=\cdots =\ED_{\mc{C}}(x_m)=0\quad\Longrightarrow \quad\ED_{\mc{C}}( x_1\cdots x_m)=0.\end{equation} A family of subsets of $\mc{B}$ is  said to be {\rm free with amalgamation over $\mc{C}$} if the subsets are contained in $*$-subalgebras which are  free with amalgamation over $\mc{C}$.
\end{Def}

%In the context of a $(\rho_1,\ldots,\rho_d)$-rectangular probability space, it is easy to see that a family  $(\chi_\alpha)_{\alpha\in A}$ of sets of simple elements are free with amalgamation over $\mc{D}$ if one has: for all $m\geq 1$, for all $\alpha_1\neq \cdots \neq \alpha_m\in A$, for all $x_1,\ldots,x_m$ elements of the $*$-algebras respectively generated by $ \chi_{\alpha_1},\ldots, \chi_{\alpha_m}$, one has  \begin{equation}\label{10.1.05.3}\ED(x_1)=\cdots =\ED(x_m)=0\quad\Longrightarrow \quad\ED( x_1\cdots x_m)=0.\end{equation}

\begin{Def} In  a $(\rho_1,\ldots,\rho_d)$-rectangular probability space,  the {\rm $\D$-distribution} of a family  $(a_j)_{j\in J}$ of simple elements is the function which maps any polynomial $P$ in the noncommutative variables $(X_j,X_j^*)_{j\in J}$, with coefficients in $\D$,  to $\ED(P(a_j,a_j^*)_{j\in J})$. \end{Def}

The last  notion we have to introduce is the {\it convergence in $\D$-distribution.} Note that the $\D$'s of different spaces can all be identified with $\C^d$.
\begin{Def} $\bullet$ If for all $n$, $(\A_n,p_{1,n},\ldots,p_{d,n},\vfi_{1,n},\ldots,\vfi_{d,n})$ is a $(\rho_{1,n},\ldots,\rho_{d,n})$-\pro space \st  $$(\rho_{1,n},\ldots,\rho_{d,n})\ninf (\rho_1,\ldots,\rho_d),$$ a family $(a_j(n))_{j\in J}$ of simple elements of $\A_n$ is said to {\rm converge in $\D$-distribution}, when $n$ goes to infinity,   to a family $(a_j)_{j\in J}$ of elements of $\A$ if the $\D$-distributions converge pointwise. 
\\ 
$\bullet$ In the context of example \ref{camille.tango.triste}, if for all $n$,  $a_j(n)$'s are $n\tii n$ random matrices, {\rm convergence in $\D$-distribution in probability} of the family  $(a_j(n))_{j\in J}$ to $(a_j)_{j\in J}$ is the convergence in probability, when $n\to\infty$, of the random variable $\ED(P(a_j(n),a_j(n)^*)_{j\in J})$  to $\ED(P(a_j,a_j^*)_{j\in J})$ for all polynomial $P$ in the noncommutative variables $(X_j,X_j^*)_{j\in J}$.\end{Def}

\subsection{Statement of the theorems about random matrices}
Let, for $s\in \n, k, l\in[d]$, $n\geq 1$, $R(s,k,l,n)$ be a $q_k(n)\tii q_l(n)$ random matrix such that for all $n$,  the family  $(R(s,k,l,n))_{s,k,l}$ is independent.  We suppose that for  all $s\in \n, k\neq l\in [d]$, $R(s,k,l,n)$ is bi-unitarily invariant and that its singular law  converge  in probability to a deterministic \pro measure. We also  suppose  that for  all $s\in \n$ odd, for all $ k\in [d]$, $R(s,k,k,n)$ is unitary uniform and that for  all $s\in \n$ even, for all $ k\in [d]$, $R(s,k,k,n)$ is hermitian, unitarily invariant and that  its spectral law converge in probability  to a deterministic \pro measure. Consider also, for $n\geq 1$,  a family $(C(i,n))_{i\in I}$ of diagonal deterministic matrices, each $C(i,n)$ having the size $q_{k_i}(n)\tii q_{l_i}(n)$, for a certain $(k_i,l_i)\in [d]\tii [d]$, which contains $p_1(n),\ldots, p_d(n)$ and is stable under product and taking adjoints. We suppose that for all $i\in I$ \st $k_i=l_i$, the normalized trace of $C(i,n)$
%$$\ff{q_{k_{i}}(n)}\Tr C(i,n)$$ 
has a finite limit when $n$ goes to infinity.

We suppose that for the operator norm $||\cdot||$ associated to canonical  hermitian norms these matrices are uniformly bounded  (in $n,s,k,l,i$).

The hypothesis of our main theorems are not the same if one supposes all $\rho_k$'s to be  positive or not. In the following theorem, all $\rho_k$'s are supposed to be positive.
\begin{Th}[Case where all $\rho_k$'s are positive]\label{14.1.05.1}
Under these hypothesis,  the family  $(\widetilde{R}(s,k,l,n)\, ;$ $s\in \n, k, l\in[d])\cup(\widetilde{C}(i,n)\ste i\in I)$ converges in $\D$-distribution in \pro to a  family   which is free with amalgamation over $\D$.
\end{Th}
%\begin{rmq} The hypothesis that $C(i,n)$'s are diagonal can be removed. Indeed, we only use it to apply case $1^\circ)$ of proposition \ref{TX1(n).merle.com}, and in the case where all $\rho_k$'s are positive, this proposition is an immediate consequence of  proposition 3.12 of \cite{fbg.free.amalg}, where constant matrices are not supposed to be diagonal. But as the author is writing this paper, \cite{fbg.free.amalg} has not been refereed yet.\end{rmq}

In the following theorem, we suppose one of the $\rho_k$'s, say $\rho_1$, to be zero. We supposed that at most one of the $\rho_k$'s is zero, so $\rho_2,\ldots,\rho_d$ are positive. 

 \begin{Th}[Case where $\rho_1=0$, other $\rho_k$'s are positive]\label{18.1.05.1} 
Assume moreover that for all $s,k\neq l,n$, $R(s,k,l,n)$ has deterministic singular values and that for all $s\in \n$ even, $k\in[d]$, $R(s,k,k,n)$ has a deterministic spectrum.   Then the family  $(\widetilde{R}(s,k,l,n)\ste s\in \n, k, l\in[d])\cup(\widetilde{C}(i,n)\ste i\in I)$ converges in $\D$-distribution in \pro to a  family   which is free with amalgamation over $\D$.
\end{Th}

\begin{cor}[Product of rectangular matrices with asymptotic null ratio]\label{16.1.08.1}Let $q_1(n), q_2(n)$ be sequences which tend to infinity in such a way that $q_1(n)/q_2(n)$ tends to zero. If, for all $n$, $M_1(n),M_2(n)$ are independent $q_1(n)\tii q_2(n)$ bi-unitarily invariant rectangular random matrices with deterministic singular values and limit singular laws, then for all $p\geq 1$, $i_1,\ldots, i_p\in \{1,2\}$, 

$\bullet$ $M_{i_1}(n)^*M_{i_2}(n)\cdots M_{i_{2p-1}}(n)^*M_{i_{2p}}(n)$ has null limit normalized trace,

$\bullet$ $M_{i_1}(n)^*M_{i_2}(n)\cdots M_{i_{2p-1}}(n)^*M_{i_{2p}}(n)$ has null limit normalized trace,
except if any matrix $M_1(n)$ or $M_2(n)$ in the product is followed by  respectively $M_1(n)^*$ or $M_2(n)^*$,  in which case the limit can be computed using the well known asymptotic freeness of square random matrices.
\end{cor}

\subsection{Comparison with already existing results}\label{adorable.cecile.12.05}

There already exists some results which allow to compute asymptotics of normalized traces of products
of rectangular random matrices: on one hand the results about Wishart matrices and on the other hand the theorem  4.1 of \cite{Sh96}.

Theorems \ref{14.1.05.1} and \ref{18.1.05.1}
allow us to compute the limit (for convergence in probability) of all normalized traces of matrices of the type \begin{equation}\label{VIH=danger.12.05}
R(s_1,k_1,l_1,n)^{\eps_1}C(i_1,n)R(s_2,k_2,l_2,n)^{\eps_2}C(i_2,n)\cdots  C(i_{p-1},n)R(s_p,k_p,l_p,n)^{\eps_p}C(i_p,n),\end{equation} with $p\geq 1,$ $s_1,\ldots, s_p\in \n$, $k_1,l_1, \ldots, k_p,l_p\in [d]$,  $\eps_1,\ldots, \eps_p\in \{1,*\}$, and $i_1,\ldots, i_{p-1}\in I$ such that the product is possible and is a square matrix.

The Wishart matrices are matrices of the type $RR^*$, where $R$ is a rectangular gaussian matrix. In   \cite{hiai}, \cite{capcas.iumj} or \cite{capdonmart}, the authors give results of asymptotic freeness for such random matrices. In the case we are interested in, where the random matrix $R$ is bi-unitarily invariant, this asymptotic freeness is only a consequence of the asymptotic freeness of unitarily invariant hermitian random matrices.
With our notations,  it  does not allow to suppose that $\rho_1=0$, and to consider products {\bf where rectangular matrices are not immediately followed by their adjoint},   as the one of (\ref{VIH=danger.12.05}).

Theorem  4.1 of the paper \cite{Sh96} by Shlyakhtenko allows, embedding rectangular matrices in larger self-adjoint matrices ($R\mapsto \begin{bmatrix}0&R\\ R^*&0\end{bmatrix}$), to consider products of the type of the one of (\ref{VIH=danger.12.05}), but {\bf only when the $\rho_i$'s are positive}, the matrices $R(s,k,l,n)$ 	are {\bf gaussian}, and the deterministic matrices are {\bf diagonal} and have, in a certain sense, {\bf a limit in $L^\infty[0,1]$}. Moreover,     this theorem states the convergence of the {\bf expectations} of the normalized traces, but not their convergence in probability (even though this could it is possible to expect to pass, using the concentration results of \cite{guoinnet+zeitouni.2000}, from convergence in expectation to convergence in probability in Shlyakhtenko's result).

As an example,  no general result had been proved about the convergence, when $n,p\to\infty$ \st $\f{n}{p}\to \la\geq 0$, of the singular law of $M(n,p)+N(n,p)$, where $M(n,p), N(n,p)$ are independent bi-unitarily invariant $n\tii p$ random matrices with limit singular laws. We did not even know if this limit existed, even though some computations have been done in some particular cases (\cite{dozier-silver04}, \cite{hln04}, \cite{hln05})). Indeed, the computation, by the moments method, of the singular law of $M(n,p)+N(n,p)$ involves the computation of the normalized trace of  $$ [(M(n,p)+N(n,p))(M(n,p)+N(n,p))^*]^k,$$the expansion of which contains products where some rectangular matrices are not followed by their adjoints.

\subsection{Proof of theorem \ref{14.1.05.1} (case where all $\rho_k$'s are positive)}

Let us introduce, on the space of $n\tii n$ matrices, the Schatten $p$-norms $||.||_p$ ($p\in [1,\infty]$), defined, for $p<\infty$,  by $$||X||_p=(\tr |X|^p)^\ff{p},$$ and $||.||_\infty=||.||$.  They satisfy classical H\"older inequalities (\cite{Nelson}). 

{\it Step I. }   Fix $k,l\in [d]$ and consider, for $n\geq 1$,   a  bi-unitarily invariant $q_k(n)\tii q_l(n)$ random matrix $R(n)$ \st the sequence $(R(n))$ is uniformly bounded for $||.||$ by $C>0$ and \st  the singular law of $R(n)$ converge in \pro to a deterministic \pro measure. Note first that, by Carleman criterion,  it implies that its singular law converges in \pro to a \pro measure $\mu$.  It is easy, transferring the proof of lemma 4.3.10 p. 160 of \cite{hiai}, to prove that the distribution of $R(n)$  can be realized  as the distribution of $U(n)\Lambda(n)V(n)$ where $U(n)$, $\Lambda(n)$, $V(n)$ are independent, $U(n)$, $V(n)$ are respectively $q_k(n)\tii q_k(n)$, $q_l(n)\tii q_l(n)$ uniform random unitary matrices, and $\Lambda(n)$ is a rectangular diagonal matrix, $$\Lambda(n)=[\delta_i^j\la(i,n)]_{\substack{1\leq i\leq q_k(n)
\\ 1\leq l\leq q_l(n)
}}\textrm{ with } 0\leq \la(1,n)\leq\cdots \leq \la(\min (q_k(n),q_l(n)),n).$$   Choose deterministic $0\leq \xi(1,n)\leq\cdots \leq \xi(\min (q_k(n),q_l(n)),n)\leq C$  \st the uniform distribution on the $\xi(1,n)$'s converges weakly to  $\mu$. Now set $$\Xi(n):=[\delta_i^j\la(i,n)]_{\substack{1\leq i\leq q_k(n)
\\ 1\leq j\leq q_l(n)
}}\textrm{ and } N(n):=U(n)\Xi(n)V(n).$$

Let us prove that for all $P$ polynomial in noncommutative random variables $X,X^*$, for all $p\in [1,\infty)$, 
$ ||P(\widetilde{R}(n))-P(\widetilde{N}(n))||_p$ converges in \pro to zero. 

By lemma 4.3.4 p. 152 of \cite{hiai}, for all $p\geq 1$,  $$  ||\widetilde{R}(n)^*-\widetilde{N}(n)^*||_p=||\widetilde{R}(n)-\widetilde{N}(n)||_p=||\widetilde{\Lambda}(n)-\widetilde{\Xi}(n)||_p$$ converges in \pro to zero. So, for all $m\geq 1$, $\eps_1,\ldots, \eps_m\in\{.,*\}$, 
\begin{eqnarray*}&||\widetilde{R}(n)^{\eps_1}\cdots \widetilde{R}(n)^{\eps_m}-\widetilde{N}(n)^{\eps_1}\cdots \widetilde{N}(n)^{\eps_m}||_p&\\
\leq &\sum_{l=1}^m||\widetilde{R}(n)^{\eps_1}\cdots \widetilde{R}(n)^{\eps_{l-1}}(\widetilde{R}(n)^{\eps_l}-\widetilde{N}(n)^{\eps_l})\widetilde{N}(n)^{\eps_{l+1}}\cdots \widetilde{N}(n)^{\eps_m}||_p&\\
\leq &\sum_{l=1}^mC^{m-1}||\widetilde{R}(n)^{\eps_l}-\widetilde{N}(n)^{\eps_l})||_p,&\end{eqnarray*} which proves the result. 

Note that in the case where $R$ is a uniformly bounded unitarily invariant $q_k(n)\tii q_k(n)$ hermitian random matrix, the moments of the spectral law of which converge in probability, then the same work can be done, replacing $V(n)$ by $U(n)^*$ and singular values by eigenvalues, and the same conclusion holds:  for all $P \in \C[X]$, for all $p\in [1,\infty)$, 
$ ||P(\widetilde{R}(n))-P(\widetilde{N}(n))||_p$ converges in \pro to zero. 
%{\it Step II. } Note that in our case, the boundedness hypothesis and the convergence in \pro of the moments of singular (or spectral) laws implies the weak convergence in \pro of  singular (or spectral) laws. Then consider, for all $s,k,l,n$, a random matrix $N(s,k,l,n)$ defined from $R(s,k,l,n)$ like $N(n)$ from $R(n)$ in the step I.  
%Once the result proved for the $R(s,k,l,n)$'s replaced by the $N(s,k,l,n)$'s, it can be deduced for $R(s,k,l,n)$'s. Indeed, 
% if $m\geq 1$, $(s_j,k_j,l_j)_{j\in [m]} \in (\n\tii [d]\tii [d])^m$, $\eps_1,\ldots, \eps_m\in\{.,*\}$, and $i_0,\ldots, i_m\in I$, 
% the normalized trace of any diagonal block of the difference between the normalized trace of $$C(i_0,n)\widetilde{R}(s_1,k_1,l_1,n)^{\eps_1}C(i_1,n)\cdots C(i_m,n)\widetilde{R}(s_m,k_m,l_m,n)^{\eps_m}C(i_m,n)$$ and  $$C(i_0,n)\widetilde{N}(s_1,k_1,l_1,n)^{\eps_1}C(i_1,n)\cdots C(i_m,n)\widetilde{N}(s_m,k_m,l_m,n)^{\eps_m}C(i_m,n)$$ has absolute value less or equal than the norm $||.||_1$ of their difference times of of the $n/\!q_k(n)$, which is less or equal, by a decompositions like in step I, than $$\ds\underset{k\in [d]}{\max }\f{n}{q_k(n)}\sum_{j=1}^mC^{2m}||\widetilde{R}(s_j,k_j,l_j,n)^{\eps_j}-\widetilde{N}(s_j,k_j,l_j,n)^{\eps_j})||_1,$$ which converges in \pro to zero, by step I and because all $\rho_k$'s are positive. So let us prove the result for $N(s,k,l,n)$'s.

{\it Step II. } 
By hypothesis,  there exists $i_1,\ldots, i_d\in I$ \st for all $n,k$, $C(i_k,n)=p_k(n)$  and there exists a map $\beta$ from $I^2$ to $I$ and an involution $*$ from $I$ \st for all $n\geq 0, i,j\in I$ \st $l_i=k_j$, \begin{equation}\label{14.1.05.3}C(i,n)C(j,n)=C(\beta(i,j),n),\quad C(i,n)^*=C(*(i),n).\end{equation}  Consider, in a  $(\rho_1,\ldots,\rho_d)$-rectangular probability space $(\A,p_1,\ldots,p_d,\vfi_1,\ldots,\vfi_d)$, a family  $a(s,k,l)$ ($s\in \n, k, l\in[d]$), $c(i)$  ($i\in I$) of elements of $\A$, with $c(i_1)=p_1$,..., $c(i_d)=p_d$, whose $\D$-distribution is defined by the following rules:
\begin{itemize}
\item[(i)] For all $s\in \n, k, l\in[d]$,  $a(s,k,l)\in p_k\A p_l$, and for all $i\in I$, $c(i)\in p_{k_i}\A p_{l_i}$.
\item[(ii)] For all $i,j\in I$ \st $l_i=k_j$,  $c(i)c(j)=c(\beta(i,j)), c(i)^*=c(*(i)).$
\item[(iii)] For all $s\in \n, k\neq l\in[d]$, for all $r\geq 0$, $$\ds\vfi_k[(a(s,k,l)a(s,k,l)^*)^r]=\underset{n\to\infty}{\lim}\tr (R(s,k,l,n) R(s,k,l,n)^*)^r,$$ $$\vfi_l[(a(s,k,l)^*a(s,k,l))^r]=\underset{n\to\infty}{\lim}\tr (R(s,k,l,n)^* R(s,k,l,n))^r.$$
\item[(iv)]  For all $s\in \n$ odd, $k\in[d]$, $a(s,k,k)$ is unitary in $p_k\A p_k$, and for all $r\in\z-\{0\}$, $\vfi_k(a(s,k,k)^r)=0$.
%is self-adjoint and for all $r\geq 0$, $$\ds\vfi_k[a(s,k,k)^r]=\underset{n\to\infty}{\lim}\tr R(s,k,l,n)^r,$$
\item[(v)]  For all $s\in \n$ even, $k\in[d]$, $a(s,k,k)$ is self-adjoint and for all $r\geq 0$, $$\ds\vfi_k[a(s,k,k)^r]=\underset{n\to\infty}{\lim}\tr R(s,k,k,n)^r,$$%the image, by $\vfi_k$, of any of its powers $a(s,k,k)^{r}$ is the limit of the normalized trace of $R(s,k,k,n)^r$. 
%\item[(v)]  For all $s\in \n$ even, $k\in[d]$, $a(s,k,k)$ is unitary in $p_k\A p_k$, and for all $r\in\z-\{0\}$, $\vfi_k(a(s,k,k)^r)=0$.
%in $a(s,k,k), a(s,k,k)^*$ is the limit of the normalized trace of the corresponding word in  $R(s,k,k,n),R(s,k,k,n)^*$.
\item[(vi)] For all $i\in I$ \st $k_i=l_i$, $\vfi_{k_i}(c(i))$ is the limit of the normalized trace of $C(i,n)$. 
\item[(vii)] The family  $\{a(s,k,l)\}$ ($s\in \n, k, l\in[d]$), $\{c(i)\ste i\in I\}$ is free with amalgamation over $\D$. 
\end{itemize}

Note that such a space and such a family exist, they are given by the free product with amalgamation (\cite{voic1},\cite{spei98}). 

 Let us denote by $\mc{C}$ the linear span of the $c(i)$'s. By (\ref{14.1.05.3}) and  (ii), for all $n$, there is a $*$-algebra morphism $\psi_n$ from $\mc{C}$ to the linear span $\mc{C}_n$ of the $\widetilde{C}(i,n)$'s. By (vi), for all $c\in \mc{C}$, for all $k\in [d]$, the normalized trace of the $k$-th diagonal block of $\psi_n(c)$ tends to $\vfi_k(c_{kk})$. 

Since $1\in \mc{C}$, it suffices to prove that for all $m\geq 0$, for all $(s_j,k_j,l_j)_{j\in [m]} \in (\n\tii [d]\tii [d])^m$, for all $P_1,\ldots,P_m$ polynomials in the noncommutative variables $X,X^*$, for all  $d(0),\ldots, d(m)\in \mc{C}$,    $\ED [\prod_{j=1}^m P_j(\widetilde{R}(s_j,k_j,l_j,n))\psi_n(d(j))]$ converges in \pro to $\ED[d(0)\prod_{j=1}^m P_j(a(s_j,k_j,l_j))d(j)].$ 

{\it Step III. } By the formula $\Tr XY=\Tr YX$ and since all $\rho_k$'s are positive,  it is easy to see that we can remove $d(0)$.

{\it Step IV. } Then, we claim that it suffices to prove it when   the $R(s,k,l,n)$'s are replaced by the $N(s,k,l,n)$'s, where the $N(s,k,l,n)$'s, are defined from the  $R(s,k,l,n)$'s like $N(n)$ from $R(n)$ in the step I.  Indeed, the normalized trace of any diagonal block of $$\ds\prod_{j=1}^m P_j(\widetilde{R}(s_j,k_j,l_j,n))\psi_n(d(j))-\ds\prod_{j=1}^m P_j(\widetilde{N}(s_j,k_j,l_j,n))\psi_n(d(j))$$
has an absolute value less or equal than the $||\cdot ||_1$ norm  of their difference times one of the $n/\!q_k(n)$, which is less or equal, by a decomposition like in step I, than $$\ds\underset{k\in [d]}{\max }\f{n}{q_k(n)}\sum_{j=1}^mC^{m}D^{m-1}||P_j(\widetilde{R}(s_j,k_j,l_j,n))-P_j(\widetilde{N}(s_j,k_j,l_j,n))||_1,$$ where $D=\underset{\substack{|z|,|z'|\leq C\\ j\in [m]}}{\max}|P_j(z,z')|$
which converges in \pro to zero, by step I and because all $\rho_k$'s are positive.

So let us prove that  for all $k\in [d]$,
 the normalized trace of the $k$-th diagonal block of $$\ds\prod_{j=1}^m P_j(\widetilde{N}(s_j,k_j,l_j,n))\psi_n(d(j))$$  converges in 
\pro to  the $k$-th coordinate of $\ED[\prod_{j=1}^m P_j(a(s_j,k_j,l_j))d(j)]. $

{\it Step V. } We will prove it by induction on $m$. If $m=0$, the result is clear. In the other case, let us suppose the result to be proved to the ranks $0,\ldots, m-1$. Let us denote $x(j)=P_j(a(s_j,k_j,l_j))$ and $X(j,n)= P_j(\widetilde{N}(s_j,k_j,l_j,n))$.  Since $\mc{C}$ is an algebra containing $p_1,\ldots,p_d$, if, for a certain $j\in[m]$, one would replace $X(j,n)$ by one of the $p_k(n)$'s and $x(j)$ by the corresponding $p_k$, then the result would follow from the induction hypothesis (and from step III if $j=1$). Thus, by linearity, one can, for all $j\in[m]$, add a linear combination of the $p_k(n)$'s to $X(j,n)$ and the same linear combination of the $p_k$'s to $x(j)$.  Therefore we can assume that for all $j$, $E(x(j))=0$. By linearity, we can moreover suppose that  for all $j$, $x(j)$ is a simple element. Note that these suppositions imply that for all $n,j$, $X(n,j)$ has only one non zero block, $\ED(X(j,n))$ converges in \pro to zero.

From now on, we use the cyclic order on $[m]$. This means that to put the index $m+1$ on an element amounts to put the index $1$. As simple elements, the $x(j)$'s have the following property: the product of any of the $x(j)$'s by one of the $p_k$'s is either zero or $x(j)$, and the same holds for $X(j,n)$'s with $p_k(n)$'s. Thus, if, for a certain $j\in[m]$ \st $(s_j,k_j,l_j)=(s_{j+1},k_{j+1},l_{j+1})$, one would replace $d(j)$ by a linear combination of the $p_k(n)$'s, then the result would follow from the induction hypothesis (and from  step III if $j=m$).  So one can suppose that  for all $j\in[m-1]$ \st $(s_j,k_j,l_j)=(s_{j+1},k_{j+1},l_{j+1})$, $E(d(j))=0$.

At last, for $j\in [m]$ \st  $(s_j,k_j,l_j)\neq (s_{j+1},k_{j+1},l_{j+1})$, one can write $$d(j)=d'(j)+\sum_{k\in[d]}\la_k(j)p_k, \textrm{ with $\ED(d'(j))=0$.}$$ Then, by linearity, one can suppose that for all such $j$,  $E(d(j))=0$ or $\exists k\in [d], d(j)=p_k.$ %But as noted above,  the product of any of the $x(j)$'s by one of the $p_k$'s is either zero or $x(j)$, and the same holds for $X(j,n)$'s with $p_k(n)$'s. So the previous assumption can be changed into $$\ED(d(j))=0\textrm{ or } d(j)=1.$$

To conclude, we have shown that we only have to prove that  $\ED[\prod_{j=1}^m X(j,n)\psi_n(d(j))]$ converges in probability to $\ED[\prod_{j=1}^m x(j)d(j)]$  under the hypothesis that all $x(j)$'s are simple elements and satisfy $\ED(x(j))=0$, and   for all $j\in [m]$,  
$$\ED(d(j))=0\textrm{ or }[d(j)\in\{p_1,\ldots,p_d\} \textrm{ and }(s_j,k_j,l_j)\neq (s_{j+1},k_{j+1},l_{j+1}) ].$$
Note that under this hypothesis, by freeness with amalgamation over $\D$,  $\ED[\prod_{j=1}^m x(j)d(j)]=0 $.

{\it Step VI. } For all $s,k,l,n$ with $k\neq l$, the matrix $N( s,k,l,n)$ has been defined by the work of step I, so we can write $$N(s,k,l,n)=U(s,k,l,n)\Xi(s,k,l,n)V(s,k,l,n),$$ with $U(s,k,l,n)$, $V(s,k,l,n)$ uniform random unitary matrices with respective sizes $q_k(n)$, $q_l(n)$, and $\Xi(s,k,l,n)$ a deterministic diagonal $q_k(n)\tii q_l(n)$ matrix. In the same way, when $s$ is even, $N(s,k,k,n)$ can be written $$N(s,k,l,n)=U(s,k,k,n)\Xi(s,k,k,n)U^*(s,k,k,n),$$ and when $s$ is odd, $N(s,k,k,n)$ is a uniform random matrix of $\mathbb{U}_k(n)$. 

Consider $j\in [m]$. If  $k_j\neq l_j$, then 
$X(j,n)$ is of one of the 4 following types: 

\begin{enumerate}\item $\widetilde{U}(s_j,k_j,l_j,n)\widetilde{\Delta}_j \widetilde{V}(s_j,k_j,l_j,n),$ where $\Delta_j$ is a $q_{k_j}(n)\tii q_{l_j}(n)$ diagonal matrix whose diagonal terms are the images of the diagonal terms of  $\Xi(s_j,k_j,l_j,n)$ by an odd polynomial.
\\
\item $\widetilde{V}(s_j,k_j,l_j,n)^*\widetilde{\Delta}_j \widetilde{U}(s_j,k_j,l_j,n)^*,$ where $\Delta_j$ is a $q_{l_j}(n)\tii q_{k_j}(n)$ diagonal matrix whose diagonal terms are the images of the diagonal terms of  $\Xi(s_j,k_j,l_j,n)$ by an odd polynomial.
\\
\item $\widetilde{U}(s_j,k_j,l_j,n)\widetilde{\Delta}_j\widetilde{ U}(s_j,k_j,l_j,n)^*,$ where $\Delta_j$ is a $q_{k_j}(n)\tii q_{k_j}(n)$ diagonal matrix whose diagonal terms are the images of the diagonal terms of  $\Xi(s_j,k_j,l_j,n)$ by an even polynomial.
\\
\item $\widetilde{V}(s_j,k_j,l_j,n)^*\widetilde{\Delta}_j\widetilde{ V}(s_j,k_j,l_j,n),$ where $\Delta_j$ is a $q_{l_j}(n)\tii q_{l_j}(n)$ diagonal matrix whose diagonal terms are the images of the diagonal terms of  $\Xi(s_j,k_j,l_j,n)$ by an even polynomial.
\end{enumerate}
And if  $k_j=l_j$, when  $s_j$ is even, only the last type is possible, and when $s_j$ is odd, $X(j,n)$ is a linear combination of positive or negative powers of $\widetilde{N}(s_j,k_j,l_j,n)$. By linearity again, this combination can be replaced by a non zero power of  $\widetilde{N}(s_j,k_j,l_j,n)$. 

Note that in all types previously presented, the constant diagonal matrix $\Delta_j$ satisfies $\ED(\widetilde{\Delta}_j)\ninf 0.$
Moreover, for all $n$, the family   of random unitary matrices  arising in these decompositions are independent,   so the result follows from case $1^\circ)$ of proposition \ref{TX1(n).merle.com} and remark \ref{TX1(n).merle.com'}.

\subsection{Proof of theorem \ref{18.1.05.1} (case where $\rho_1=0$, all other $\rho_k$'s are positive)} 
We shall follow closely the proof of theorem \ref{14.1.05.1} (case where all $\rho_k$'s are positive) given in the previous section, but because of the weaker hypothesis (singular values of non hermitian random matrices and spectrum of hermitian ones are now supposed to be deterministic), we shall skip steps I and IV, in which non hermitian (resp. hermitian) random matrices where approximated by random matrices with deterministic singular values  (resp. deterministic spectrum). However, the hypothesis    $\rho_1=0$ will make things slightly harder.  

Copying step II of the previous proof, we have to prove, for all $m\geq 0$,  the following proposition: $H(m)$:= "for all $(s_j,k_j,l_j)_{j\in [m]} \in (\n\tii [d]\tii [d])^m$,  for all $P_1,\ldots,P_m$ polynomials in the noncommutative variables $X,X^*$,
for all $d(0),\ldots, d(m)\in \mc{C}$, for all $k\in [d]$, $$\ED[\psi_n(d(0))\prod_{j=1}^m P_j(\widetilde{R}(s_j,k_j,l_j,n))\psi_n(d(j))]$$
converges in \pro to $\ED[d(0)\prod_{j=1}^m P_j(a(s_j,k_j,l_j))d(j)].$" 

Note first that by linearity, it suffices to prove it when $d(0), d(m)$ are supposed to be simple, which implies that there exists $u,v\in [d]$ \st 
$p_ud(0)=d(0)$, $d(m)p_v=d(m).$
Let us denote by $H(m,u,v)$ the proposition $H$ in the particular case where $p_ud(0)=d(0)$, $d(m)p_v=d(m).$ Let us also denote by $H'(m)$ the proposition $H(m)$ in the particular case where $d(0)=1$.

Note first that $H(m,u,v)$ is obvious when $u\neq v$. Moreover, like in step III of the previous section, we see that for all $u\in [d]-\{1\}$, $H'(m)$ implies $H(m,u,u)$. Thus we only have to prove that for all $m\geq 0$, $H'(m)$ and $H(m,1,1)$ hold. 

Let us prove it by induction on $m$. For $m=0$, the result is immediate. Consider $m\geq 1$ \st $H'(l)$ and $H(l,1,1)$ hold for all $l<m$. 

Let us first prove $H'(m)$. Consider  $(s_j,k_j,l_j)_{j\in [m]} \in (\n\tii [d]\tii [d])^m$,  $P_1,\ldots,P_m$ polynomials in the noncommutative variables $X,X^*$,
$d(1),\ldots, d(m)\in \mc{C}$.  Let us prove that $$\ED[\prod_{j=1}^m P_j(\widetilde{R}(s_j,k_j,l_j,n))\psi_n(d(j))]$$  converges in \pro to $\ED[\prod_{j=1}^m P_j(a(s_j,k_j,l_j))d(j)]$.  Let us denote $x(j)=P_j(a(s_j,k_j,l_j))$ and $X(j,n)= P_j(\widetilde{N}(s_j,k_j,l_j,n))$.  
Similarly to step V of the previous section, by the induction hypothesis and the fact that  $H'(m-1)$ implies $H(m-1,u,u)$, we can suppose that for all $n,j$, $X(n,j)$ has only one non zero block and 
 $\ED(X(j,n))$ converges in probability to zero. 
We can also suppose that for all $j\in [m]$ such that  with the cyclic order, $(s_j,k_j,l_j)=(s_{j+1},k_{j+1},l_{j+1})$,  $E(d(j))=0$  or there exists $k\in [d]$ \st $d(j)=p_k.$ 

Under these assumptions,  we only have to prove that  $\ED[\prod_{j=1}^m X(j,n)\psi_n(d(j))]$  converges in 
\pro to  $\ED[\prod_{j=1}^m x(j)d(j)]=0.$ Following step VI of the previous section, it appears to be an application of 
 case $1^\circ)$ of proposition \ref{TX1(n).merle.com} and remark \ref{TX1(n).merle.com'}.

Similarly,  $H(m,1,1)$ is an application of case $2^\circ)$ of proposition \ref{TX1(n).merle.com}.

\section{The rectangular free additive convolution}\label{ma.raison.vient.tjrs.me.dire}

For $\nu$ probability measure on $\R$,  denote by $\tilde{\nu}$ the {\it symmetrization} of $\nu$, which is the  \pro measure defined by  $\tilde{\nu}(B)=\ff{2}(\nu(B)+\nu(-B))$ for all Borel set $B$.

One of the interests of modeling the asymptotic behavior of square random matrices with free probabilities is the possibility to compute the asymptotic spectral law of an hermitian matrix which is a function of several independent random matrices. In particular (\cite{voic91}), if $M(1,n)$, $M(2,n)$ are independent hermitian unitarily invariant random matrices whose spectral laws tend respectively to $\nu_1,\nu_2$ when their dimension $n$ goes to infinity, then the spectral law of $M(1,n)+M(2,n)$  tends to $\nu_1\bxp \nu_2$ (free additive convolution of $\nu_1,\nu_2$), and the spectral law of $\sqrt{M(1,n)}M(2,n)\sqrt{M(1,n)}$ (when the matrices are positive) tends to   $\nu_1\bxt \nu_2$ (free multiplicative convolution of $\nu_1,\nu_2$). In the same way (combine Theorem 4.3.11 of \cite{hiai} and Propositions 3.5, 3.6 of \cite{haag2}), if     the matrices $M(1,n),M(2,n)$ are square,  bi-unitarily invariant (and still independent), with asymptotic singular distributions $\nu_1$, $\nu_2$,  then the symmetrization of the  singular distribution of  $M(1,n)+M(2,n)$ tends to  $\tilde{\nu}_1\bxp \tilde{\nu}_2$, and the push-forward, by the function $x\to x^2$, of the singular distribution of $M(1,n)M(2,n)$ is the free multiplicative convolution of the push-forwards, by $x\to x^2$, of $\nu_1,\nu_2$. 

Now, suppose that  $M(1,n)$, $M(2,n)$, instead of being square, are rectangular (independent and bi-unitarily invariant) $q_1(n)\tii q_2(n)$ random matrices whose singular laws  tend to $\nu_1,\nu_2$. We keep the notations introduced in section \ref{23.1.05.132}, but suppose that $d=2$: $$q_1(n)+q_2(n)=n,\quad q_1(n), q_2(n)\ninf \infty,\quad \f{q_1(n)}{n}\ninf \rho_1\geq 0,\quad\f{q_2(n)}{n}\ninf \rho_2>0. $$  
Since for $1\leq q\leq q'$, the singular law $\nu$ of a $q\tii q'$ matrix $M$ is related to the singular law $\nu'$ of $M^*$ by $\nu'=(q'-q)/q'\delta_0+q/q'\nu,$ we can, without restriction, suppose that $q_1(n)\leq q_2(n)$.

First, the singular distribution of $M(1,n)M(2,n)^*$ tends to a distribution that can be computed by free probability theory. Indeed, for all $r\geq 1$ even, the $r$-th moment of the singular distribution of $M(1,n)M(2,n)^*$ is 
\begin{equation}\label{24.1.05.1}\f{\Tr (M(1,n)M(2,n)^*M(2,n)M(1,n)^*)^{r/\! 2}}{q_1(n)}=\f{q_2(n)\Tr (M(2,n)^*M(2,n)M(1,n)^*M(1,n))^{r/\! 2}}{q_1(n)q_2(n)},\end{equation}
which tends to $\rho_2/\rho_1$ (when $\rho_1>0$) times the $r/\!2$-th moment of $\tau_1\bxt \tau_2$, where   $\tau_1,\tau_2$
are the limit spectral distributions of $M(2,n)^*M(2,n),M(1,n)^*M(1,n)$. $\tau_1, \tau_2$ can easily be computed from $\nu_1, \nu_2$. Thus, if $\rho_1>0$, the limit of the spectral distribution of  $M(1,n)M(2,n)^*$ can be easily computed using tools of free probability. If $\rho _1=0$,  then the explicit computation of the expectation and of the variance of (\ref{24.1.05.1}) for $r=2$ gives us the convergence of the singular law of  $M(1,n)M(2,n)^*$ to the Dirac measure in zero. 
 
But until now, it wasn't possible to compute the singular  law of $M(1,n)+M(2,n)$ with the tools of free \pro theory. Indeed, for $r\geq 1$ even, its $r$-th moment is given by $$\ff{q_1(n)}\Tr ((M(1,n)+M(2,n))(M(1,n)^*+M(1,n)^*))^{r/\! 2},$$
which can be expanded into the sum,  over all $i\in \{1,2\}^r$, of \begin{equation}\label{24.1.05.3}\ff{q_1(n)}\Tr M(i_1,n)M(i_2,n)^*\cdots M(i_{r-1},n)M(i_r,n)^*.\end{equation}  Asymptotics of normalized traces like in (\ref{24.1.05.3}) cannot be computed with free probability theory, because our random matrices are not square, and we cannot reduce the problem, as for $M(1,n)M(2,n)^*$, to a problem which involves only independent square random matrices.  Thus we have to use theorem \ref{14.1.05.1} (or corollary  \ref{16.1.08.1} if $\rho_1=0$) to compute the asymptotic singular law of $M(1,n)+M(2,n)$. 

\begin{propdef}\label{15.7.4.1} Let $\mu_1,\mu_2$ be two compactly supported symmetric \pro measures on the real line. \\
(a) Let, for all $n\geq 1$, $M(1,n)$, $M(2,n)$ be independent bi-unitarily invariant $q_1(n)\tii q_2(n)$ random matrices, with deterministic singular values if $\rho_1=0$, uniformly bounded for $||.||$, and \st for all $i=1,2$, the symmetrization of the singular law of $M(i,n)$ converges in \pro to $\mu_i$. Then the symmetrization of the singular law of $M(1,n)+M(2,n)$ converges in \pro to a symmetric \pro measure on the real line, denoted  $\mu_1\arc \mu_2$, which depends only on $\mu_1,\mu_2$, and $\la:=\lim_{n\to +\infty} q_1(n)/\! q_2(n)$.  
\\
(b) $\mu_1\arc \mu_2$ is the unique  symmetric measure $\mu$ \st for all $r\in \n$ even, the $r$-th moment of $\mu$ is $\vfi_1(((a_1+a_2)(a_1+a_2)^*)^{r/\! 2})$, where $a_1$, $a_2$ are free with amalgamation over $\D$ elements of $p_1\A p_2$, $(\A,p_1,p_2,\vfi_1,\vfi_2)$ is a $(\rho_1,\rho_2)$-\pro space, and for all $i=1,2$, for all $r\in \n$ even, $\vfi_1((a_ia_i^*)^{r/\! 2})$ is the $r$-th moment of $\mu_i$.
\\
(c) The support of $\mu_1 \arc\mu_2$ is contained in the sum of the convex hulls of the supports of $\mu_1$ and $\mu_2$. 
\end{propdef}
The binary operation $\arc$ on the set of compactly supported symmetric \pro measures, called {\it rectangular free convolution with ratio $\la$},  will be extended in section \ref{perez.la.praline} to the set of symmetric \pro measures, and the same result will stay true without any hypothesis of boundedness.

\begin{pr} 
By theorem \ref{14.1.05.1} (or theorem \ref{18.1.05.1} if $\rho_1=0$), if  $a_1$, $a_2$ are elements of a $(\rho_1,\rho_2)$-\pro space $(\A,p_1,p_2,\vfi_1,\vfi_2)$ as in (b), then  $(\widetilde{M}(1,n), \widetilde{M}(2,n))$ converges in $\D$-distribution in \pro to $(a_1,a_2)$. It implies that for all $k\geq 1$, the normalized trace of $[(M(1,n)+M(2,n))(M(1,n)+M(2,n))^*]^k$ converges in \pro to a  finite and deterministic limit $m_k$ as $n$ goes to infinity. 

Let us first prove that the $m_k$'s are the moments of a \pro measure with support on $[0,(s_1+s_2)^2]$. It suffices to prove that for any real polynomial $P=\sum_kc_kX^k$ which is nonnegative on this interval,  $\sum_kc_km^k\geq 0$.   Note that since the $m_k$'s only depend of the $\D$-distribution of $(a_1,a_2)$  one can suppose that  for all $n$, for all $i=1,2$, $$M(i,n)=\ds U(i,n)[\delta_k^l\la_i(k,n)]_{\substack{
1\leq k\leq q_1(n)\\ 1\leq l\leq q_2(n)
}}V(i,n),$$
 where 
$\la_i(1,n),\ldots , \la_i(q_1(n),n)$ are  real numbers of the support of $\mu_i$  \st the uniform distribution on the $\la_i(k,n)$'s ($k=1,\ldots, q_1(n) $) converges weakly to $\mu_i$ and $U(i,n)$, $V(i,n)$, are independent respectively $q_1(n)\tii q_1(n)$, $q_2(n)\tii q_2(n)$ uniform unitary random matrices. In this case, $\sum_kc_km^k$ is the limit of the normalized trace of $$\sum_k c_k[(M(1,n)+M(2,n))(M(1,n)+M(2,n))^*]^k.$$ By hypothesis on $P$, this matrix is hermitian and nonnegative for all $n$, hence $\sum_kc_km^k\geq 0$.

Since on the set of \pro measures on a compact interval, convergence of moments is equivalent to weak convergence, (a) and (c) are proved. Moreover, (b) follows from the first paragraph of the proof.
\end{pr}

\begin{rmq}When $\la=1$,  the rectangular free convolution with ratio $\la$  is the well known additive free  convolution.
 To see it, apply the previous proposition to sequences of  square matrices. We know (Theorem 4.3.11 of \cite{hiai}) that independent bi-unitarily invariant square matrices are asymptotically free $R$-diagonal elements and  that the symmetrization of the distribution of the absolute value of the sum of two free $R$-diagonal elements is the free convolution of the symmetrizations of  the distributions of their absolute values (Proposition 3.5 of \cite{haag2}), which proves that $\arc=\bxp$ for $\la=1$.
We will see later that the free convolution with ratio $0$ is also related to the free convolution.\end{rmq}

\section{The rectangular $R$-transform}\label{cum.20.09.06}

\subsection{Relation moments-cumulants in a $(\rho_1,\rho_2)$ \pro space}

For an introduction to the general theory of cumulants in a $\D$-probability space and the particular case of a $(\rho_1,\rho_2)$)\pro space, we refer to section 2 of \cite{fbg.free.amalg}. Basically, if $(\A,p_1,p_2,\vfi_1,\vfi_2)$ is a $(\rho_1,\rho_2)$-\pro space, then for all $n$, there are two linear functions   $c_n^{(1)}, c_n^{(2)}$ from the $n$-th tensor product with amalgamation over $\D$ of $\A$ to $\C$ such that:

(a) for all $a\in p_1\A p_2$, for all $n\geq 1$, \begin{eqnarray}\vfi_1((aa^*)^n)&=&\ds\sum_{\pi\in \NCp(2n)} \prod_{\substack{V \in\pi\\ \min V \textrm{ odd}}} c_{|V|}^{(1)}(a\otimes \cdots \otimes a^*)\prod_{\substack{V \in\pi\\ \min V \textrm{ even}}} c_{|V|}^{(2)}(a^*\otimes \cdots \otimes a),\label{17.01.08.1}\\
\vfi_2((a^*a)^n)&=&\ds\sum_{\pi\in \NCp(2n)} \prod_{\substack{V \in\pi\\ \min V \textrm{ odd}}} c_{|V|}^{(2)}(a^*\otimes \cdots \otimes a)\prod_{\substack{V \in\pi\\ \min V \textrm{ even}}} c_{|V|}^{(1)}(a\otimes \cdots \otimes a^*),\label{17.01.08.3}
\end{eqnarray} where $\NCp(2n)$ is the set of noncrossing partitions of $[2n]$ where all blocks have an even cardinality. This formula can easily be deduced from equations (2.2), (2,6) and (b) of section 2.2 of \cite{fbg.free.amalg},

(b) $a\in p_1\A p_2$, for all $n$ even, we have $\rho_1c_n^{(1)}(a\otimes a^* \otimes \cdots \otimes a^*)=\rho_2c_n^{(2)}( a^*\otimes a \otimes \cdots \otimes a)$,

(c) for all $a,b\in p_1\A p_2$ free with amalgamation over $\D$, for all for all $n$ even, we have \begin{equation}\label{17.01.08.2}c_n^{(1)}((a+b)\otimes (a+b)^* \otimes \cdots \otimes (a+b)^*)=c_n^{(1)}( a\otimes a^* \otimes \cdots \otimes a^*)+c_n^{(1)}( b\otimes b^* \otimes \cdots \otimes b^*).\end{equation}

From  this, we deduce, for $\la=\rho_1/\rho_2$: 
\begin{propo}\label{micorondes.25.1.05}For $a\in p_1\A p_2$, $n$ positive integer, \begin{equation}\label{25.1.05.1}\vfi_1((aa^*)^n)=\ds\sum_{\pi\in \NCp(2n)}\la^{e(\pi)}\prod_{V\in \pi}c_{|V|}(a),\end{equation}where for all $r$ even, $c_r(a)$ denotes $c_r^{(1)}(a\otimes \cdots \otimes a^*)$, and for all $\pi$, $e(\pi)$ denotes the number of blocks of $\pi$ with even minimum.\end{propo}

\subsection{The case $\la=0$, characterization of $\arco$}\label{24.05.04.1} 
\begin{lem}\label{17.01.08.777}If $\la=0$, then for all $a\in p_1\A p_2$, for all $n\geq 1$, $c_{2n}(a)$ is the $n$-th free cumulant of the element $aa^*$ of $(p_1\A p_1, \vfi_1)$.
\end{lem}

\begin{pr} It suffices to prove that for all $n\geq 1$,  \begin{equation}\label{25.1.05.1.aprem}\vfi_1((aa^*)^n)=\ds\sum_{\pi\in \NC (n)}\prod_{V\in \pi}c_{2|V|}(a).\end{equation} By \eqref{25.1.05.1}, the left hand term is equal to $\sum_{\pi\in \NCp(2n), e(\pi)=0} \prod_{V\in \pi}c_{|V|}(a)$. Let us introduce the set $\NCd (n)$ of non-crossing partitions $\pi$ of $[2n]$ \st for all $k$ even in $[2n]$, $k-1\stackrel{\pi}{\sim}k$. Note that $\NCd (2n)$ is contained in the set of elements of $\NCp(2n)$   in which no block has even minimum. We claim that the inverse inclusion is true. Indeed, suppose the existence  of $\pi\in \NCp(2n)\backslash \NCd(2n)$ in which no block  has even minimum. Define $k=\min\{k\in[2n]\ste $\textrm{$k$ even and }$k-1\stackrel{\pi}{\nsim}k\}.$  $k$ cannot be the minimum of its block in $\pi$, and by minimality of $k$, its preceding element $j$ in this block has to be even.  Then, as $\pi$ is non-crossing, the set $\{j+1,$..., $k-1\}$ is a union of classes of $\pi$, but its cardinality is odd, which is in contradiction with $\pi\in \NCp(2n)$. 

Now, \eqref{25.1.05.1.aprem} follows easily from the fact that $\NCd(2n)$ is in correspondence with  $\NC(n)$ by the order-preserving bijection from $\NC(n)$ onto $\NCd (2n)$, that maps any noncrossing partition $\pi$ to the partitionthat links two elements $i,j$ of $[2n]$ \ssi the upper integer parts of $i/2$ and $j/2$ are linked by $\pi$. 
\end{pr}

Now we are able to give the link between the free convolution with null ratio $\arco$ and the free convolution $\bxp$:
\begin{propo}\label{18.09.06.1}The free convolution with null ratio of two compactly supported symmetric \pro measures is the unique symmetric \pro whose   push-forward by the square function is the free convolution of the push-forwards by the square function of the two probabilities.  
\end{propo}

\begin{pr} Denote, for $\nu$ \pro measure, $m_n(\nu)$ its $n$-th moment, and $\kk_n(\nu)$ its $n$-th free cumulant. Recall (\cite{spei}) that free cumulants of \pro measures  are defined by the formula \begin{equation}\label{24.05.04.3}\ds \forall n\geq 1, m_n(\nu)=\sum_{\pi\in\NC(n)}\prod_{\substack{B\in \pi}}\kk_{|B|}(\nu),\end{equation} and that the free convolution of $\nu_1,\nu_2$ compactly supported \pro measures is the only \pro with free cumulants $\kk_n(\nu_1)+\kk_n(\nu_2)$ ($n\geq 1$).  
\\
Consider two symmetric compactly supported \pro measures $\mu_1,\mu_2$. Denote by $\nu_1,\nu_2$ their respective push-forwards by the square function. It suffices to prove that for all $n\geq 1$, $m_{2n}(\mu_1\arco\mu_2)=m_{n}(\nu_1\bxp\nu_2)$, i.e. that  for all $n\geq 1$, $$\ds m_{2n}(\mu_1\arco\mu_2)=\sum_{\pi\in\NC(n)}\prod_{\substack{B\in \pi}}(\kk_{|B|}(\nu_1)+\kk_{|B|}(\nu_2)).$$
Consider  $a_1,a_2\in p_1\A p_2$ free with amalgamation over $\D$ \st $$\forall i=1,2, \forall n\geq 1,  \vfi_1((a_ia_i^*)^{n})=m_{2n}(\mu_i).$$
Then one has, for all $n\geq 1$,  \begin{eqnarray*} m_{2n}(\mu_1\arco\mu_2)&=&\vfi_1[((a_1+a_2)(a_1^*+a_2^*))^{n}]\\
&=&\sum_{\pi\in\NCd (2n)}\prod_{\substack{B\in \pi}}\underbrace{c_{|B|}(a_1+a_2)}_{=c_{|B|}(a_1)+c_{|B|}(a_2)\textrm{ by (\ref{17.01.08.2})}
}
\textrm{by  (\ref{25.1.05.1})} \end{eqnarray*}
But by lemma \ref{17.01.08.777}, for all $i=1,2$, for all $m$ even, $c_m(a_i)$ is the $m/\!2$-th free cumulant of the element $a_ia_i^*$ of $(\A_1,\vfi_1)$, i.e. the $m/\!2$-th free cumulant of $\nu_i$. 
%&=&\sum_{\pi\in\NC (n)}\prod_{\substack{B\in \pi}}(\kk_{|B|}(a_1a_1^*,\ldots,a_1a_1^*)+\kk_{|B|}(a_2a_2^*,\ldots))\\
So  \begin{eqnarray*} m_{2n}(\mu_1\arco\mu_2)&=&\sum_{\sigma\in\NC (n)}\prod_{\substack{V\in \sigma}}(\kk_{|V|}(\nu_1)+\kk_{|V|}(\nu_2))=m_{n}(\nu_1\bxp\nu_2),
\end{eqnarray*}which proves the result.
\end{pr}

\subsection{Generating series of the cumulants of an element of $p_1\A p_2$}\label{27.04.04.4} We are going to derive from (\ref{25.1.05.1}) a formula which links generating series of the sequence $\vfi_1((aa^*)^n)$ and $c_{2n}(a)$. Note that such a formula could be derived from the very general Theorem 2.2.3 of \cite{spei98}, but the work needed to translate this result to our context and tu put it into the form we shall use it is as long as this section. 

In this section, we work in the field $\C((X))$ of fractions of the ring $\C[[X]]$ of formal power series in $X$ with complex coefficients, endowed with the classical addition and product. Note that for any element $S$ of $\C[[X]]$ with null constant term,  $F\mapsto F\circ S$ is a well defined operation on $\C((X))$. If moreover, the coefficient of $X$ in $S$ is non null, then $S$ has a left and right inverse for $\circ$, denoted by $S^{<-1>}$. 
\begin{lem}\label{chouineur}  Consider a sequence $(c_{2n})_{n\in \n^*}$ of complex numbers. Define the sequence  $(m_{2n}^{(e)})_{n\in \n}$ %$(m_{2n}^{(o)})_{n\in \n}$ 
by $m_0^{(e)}=1$ %,  $m_0^{(o)}=1$, 
and for each $n\in \n^*$, $$m_{2n}^{(e)}=\sum_{\pi\in \NCp(2n)}\la^{e(\pi)}\prod_{B\in\pi}c_{|B|}.$$%,\qquad m_{2n}^{(o)}=\sum_{\pi\in NC'(2n)}\la^{o(\pi)}\prod_{B\in\pi}c_{|B|}.$$ 
Define the formal power series $$C(X)=\sum_{n\geq 1} c_{2n}X^n,\quad\quad M^{(e)}(X)=\sum_{n\geq 1}m_{2n}^{(e)}X^n.$$%,\; M^{(o)}(X)=\sum_{n\geq 1}m_{2n}^{(o)}X^n.$$
Then we have $\ds C[X  (\la M^{(e)}(X) +1)(M^{(e)}(X)+1)]=M^{(e)}(X)$. 
\end{lem}

\begin{pr} {\it Step I. } Define, for $n$ positive integer and $\pi\in \NCp(2n)$, $o(\pi)$ to be the number of blocks of $\pi$ with odd minimum. Then we have $$\la  m_{2n}^{(e)} =\ds\sum_{\pi\in \NCp(2n)} \la^{o(\pi)}\prod_{B\in\pi}c_{|B|}.$$ Indeed, if $c$ denotes the cycle $(2n\to 2n-1\to \cdots \to 2\to 1)$ of $[2n]$, then $c$ induces a permutation of $\NCp(2n)$ by $$R :\pi\in \NCp(2n)\mapsto R(\pi):=\{c(B)\ste B\in \pi\}.$$
For example, for $\pi=\{\{1, 4,5,6\},\{2,3\}\}$, $R(\pi)=\{\{ 3,4,5,6\},\{1,2\}\}$.
For $\pi\in \NCp(2n)$, for $B$ block of $\pi$, $\min c(B)$ is $(\min B)-1$ if $\min B>1$, and the second element of $B$ minus $1$ in the other case. Note that in the first case, $\min B$ and $\min c(B)$ have different parities, whereas in the second one, they have the same parities, by definition of $\NCp(2n)$. Hence we have $o(R(\pi))=e(\pi)+1$. Thus  $$\la  m_{2n}^{(e)}=\sum_{\pi\in \NCp(2n)}\la^{e(\pi)+1}\prod_{B\in\pi}c_{|B|} =\ds\sum_{\pi\in \NCp(2n)} \la^{o(R(\pi))}\prod_{B\in\pi}c_{|B|}=\sum_{\pi\in \NCp(2n)}\la^{o(\pi)}\prod_{B\in\pi}c_{|B|}.$$ 
 
{\it Step II. } Thus we have to prove that  $\ds C[X  (M^{(o)}(X) +1)(M^{(e)}(X)+1)]=M^{(e)}(X)$, where $M^{(o)}(X)=\sum_{n\geq 1}m_{2n}^{(o)}X^n,$ and $m^{(o)}_{2n}=\ds\sum_{\pi\in \NCp(2n)} \la^{o(\pi)}\prod_{B\in\pi}c_{|B|}.$

Define,  for each positive $n$ and each $j\in [n]$, the ``decomposition map'' $D_j$, from the set of partitions of $\NCp(2n)$ in which the class of $1$ has cardinality $2j$ to $\underset{\substack{l\in \n^{2j}\\ l_1+\cdots +l_{2j}=n-j}}{\cup}\NCp(2l_1)\times\cdots\times \NCp(2l_{2j})$ (the set $\NCp(0)$ is considered as a singleton, on which $e$ and $o$ have value $0$),  which maps $\pi\in NC'(n) $ to 
$(\pi_1,\ldots,\pi_{2j})$ defined in the following way:  if $\{1=k_1<\ldots<k_{2j}\}$ is the class of $1$ in $\pi$, then for each $r$, let $\pi_r$ be the restriction of $\pi$ to the interval $]k_r,k_{r+1}[$, with $k_{2j+1}=2n+1$. % $l_j=n-j-(l_1+\cdots +l_{j-1})$, $l_0=0$, and for each $r$, $\pi_r=\pi_{|]r+l_1+\cdots +l_r-1,r+l_1+\cdots +l_r ]}$.  
It is easy to see that the map $D_j$ is  well defined (recall that for all $r$, $k_r$ and $k_{r+1}$ have opposite parity), and that it is a bijection. Moreover, we have \begin{eqnarray*}e(\pi)&=&o(\pi_1)+e(\pi_2)+\cdots+o(\pi_{2j-1})+e(\pi_{2j}),\\ o(\pi)&=&1+e(\pi_1)+o(\pi_2)+\cdots+e(\pi_{2j-1})+o(\pi_{2j}).
\end{eqnarray*}
We denote by $\langle X^n\rangle P$ the coefficient of $X^n$ in a formal power series $P$ in $X$. We have to show that for each $n\geq 1$, $$\langle X^n\rangle C[X (M^{(o)}(X) +1)(M^{(e)}(X)+1)]=m_{2n}^{(e)}.$$
We have 
\begin{eqnarray*}&
\langle X^n\rangle C[X  (M^{(o)}(X) +1)(M^{(e)}(X)+1)]&\\ =&\ds\sum_{j=1}^n c_{2j}\langle X^n\rangle [X^j   (M^{(o)}(X) +1)^j(M^{(e)}(X)+1)^j]&\\ 
=&\ds\sum_{
j=1
}^n
c_{2j}\langle X^{n-j}\rangle [(M^{(o)}(X) +1)^j(M^{(e)}(X)+1)^j]&\\
=&\ds\sum_{j=1}^n c_{2j}\underset{\substack{l\in \n^{2j}\\ l_1+\cdots +l_{2j}=n-j}}{\sum}m_{2l_1}^{(o)}m_{2l_2}^{(e)}\cdots m_{2l_{2j-1}}^{(o)}m_{2l_{2j}}^{(e)}&\\
=&\ds\sum_{j=1}^n c_{2j}\!\!\!\!\underset{\substack{l\in \n^
{2j}\\ l_1+\cdots +l_{2j}=n-j}}{\sum}\underset{\substack{\pi_1\in \NCp(2l_1)\\ \vdots\\ \pi_{2j}\in \NCp(2l_{2j})}}{\sum}\la^{o(\pi_1)+e(\pi_2)+\cdots+e(\pi_{2j})}\prod_{r=1}^{2j}\prod_{B\in \pi_{r}}c_{|B|}&
\end{eqnarray*}
Hence the preliminaries about the bijections $D_j$ tell us that $$
\langle X^n\rangle C[X  (M^{(o)}(X) +1)(M^{(e)}(X)+1)]\,=\,\sum_{\pi\in \NCp(2n)}\la^{e(\pi)}\prod_{B\in \pi}c_{|B|}\,\;=\;\, m_{2n}^{(e)}.$$\end{pr}

Let us introduce two power series that will play a role in the computation of the generating series of the cumulants:\begin{eqnarray*}T(X)&=&(\la X+1)(X+1),\\
U(X)\;=\; (T-1)^{<-1>}
&=&\sum_{n\geq 1}\f{(-1)^{n-1}\la^{n-1}{2n\choose n}}{2(\la +1)^{2n-1}(2n-1)}X^n.\end{eqnarray*} $U(z)$ is the power expansion of 
$\f{-\la-1+\lf[(\la+1)^2+4\la z\ri]^{1/2}}{2\la}$ ($:=z$ if $\la=0$), 
where $z\mapsto z^{1/2}$ is the analytic version of the square root on the complement of the real non positive half line \st $1^{1/2}=1$.    

\begin{Th}\label{21.01.08.1}Consider $a\in p_1\A p_2$ and define the formal power series
$M(X)=\sum_{n\geq 1}m_{2n}(a)X^n$ and $C(X)=\sum_{n\geq 1}c_{2n}(a)X^n$. Then we have $$C=U \lf(
\f{
X
}{
\lf(
X{\scriptstyle{\times}} (T\circ M)
\ri)^{<-1>}
}
-1
\ri).$$
\end{Th}

\begin{pr} By proposition \ref{micorondes.25.1.05} and the  previous lemma, we have  $C\circ\lf(
X{\scriptstyle{\times}} (T\circ M)
\ri)=M.$ Hence $T\circ C\circ\lf(
X{\scriptstyle{\times}} (T\circ M)
\ri) =T\circ M$. Dividing by $T\circ M$ on both sides, one gets $
\f{T\circ C}{X}\circ\lf(
X{\scriptstyle{\times}} (T\circ M)
\ri)=\ff{X}.$ Then by inversion, $ \f{X}{T\circ C}\circ\lf(X{\scriptstyle{\times}} (T\circ M)\ri)=X$. At last, $
\f{X}{T\circ C}=\lf(X{\scriptstyle{\times}} (T\circ M)\ri)^{<-1>}$, which implies $ T\circ C=\f{X
}{
\lf(
X{\scriptstyle{\times}} (T\circ M)
\ri)^{<-1>}
}$, and then  $
C=U \lf(
\f{
X
}{
\lf(
X{\scriptstyle{\times}} (T\circ M)
\ri)^{<-1>}
}
-1
\ri).$
\end{pr}

%Until the beginning  of section \ref{perez.la.praline}, $\la$ will denote a number of $(0,1]$. 
\subsection{Rectangular $R$-transform of compactly supported symmetric measures}\label{supremepatience} From now on, we will denote by $(\cdot)^{1/2}$ (resp. $\sqrt{\cdot}$) the analytic version of the square root on the complementary of  the real non positive (resp. nonnegative) half line \st $1^{1/2}=1$ (resp. $\sqrt{-1}=i$).    

Consider a symmetric compactly supported \pro measure $\mu$. For $n\geq 0$, denote by $m_n(\mu)$ its $n$-th moment and by $c_{2n}(\mu)$ the number $ c_{2n}(a)$ with $a\in p_1\A p_2$, where $\A$ is a $(\rho_1,\rho_2)$-\pro space, \st for all $k\geq 1$, $\vfi_{1}((aa^*)^k)= m_{2k}(a)$.
%\begin{rmq}\label{24.05.04.2} In the case $\la=0$, by the results of subsection \ref{24.05.04.1}, it is easy to see that $c_{2n}(\mu)$ is the $n$-th free cumulant $\kk_n(\rho)$ of the push-forward $\rho$ of $\mu$ by the square function.\end{rmq}

Note that the generating series of the moments $M_\mu(X)=\sum_{n\geq 1}m_{2n}(\mu)X^n$ of $\mu$ is the power expansion of $$M_\mu(z)=\ff{\sqrt{z}} G_\mu(\ff{\sqrt{z}} )-1 \quad \quad (z\in \C\bck [0,+\infty)),$$ where $G_\mu$ is the  {\textit{Cauchy transform}} of $\mu$ on the upper half plane: $G_\mu(z)=\int_{t\in\R}\f{\ud \mu(t)}{z-t}. $ Hence $X{\scriptstyle{\times}} (T\circ M)$  is the power expansion of $$H_\mu(z):=zT\circ M_\mu(z)=\la G_\mu\lf(\ff{\sqrt{z}}\ri)^2+(1-\la)\sqrt{z}G_\mu\lf(\ff{\sqrt{z}}\ri) \quad \quad (z\in \C\bck [0,+\infty)), $$that we shall call  the {\textit{rectangular Cauchy transform with ratio $\la$}}  of $\mu$.   
Since $\mu$ is compactly supported, from the power expansion of the Cauchy transform $G_\mu$, we know that $H_\mu$ is analytic in a neighborhood of zero, that $H_\mu(0)=0$ and that $H_\mu'(0)=1$. 
It implies that we can invert the function $H_\mu$ in a neighborhood of zero  and that the inverse function $H_\mu^{-1}$ has a power series expansion  given by $\lf(X\tii (T\circ M_\mu) \ri)^{<-1>}$. Note that   $\f{z}{H_\mu^{-1}(z)}-1$ is analytic in a neighborhood of zero, with value zero in zero. At last, by theorem \ref{21.01.08.1}, $C_\mu(X):=\sum_{n\geq 1}c_{2n}(a)X^n$ is the power expansion of $$C_\mu(z):= U\lf(\f{z}{H_\mu^{-1}(z)}-1\ri)\quad\textrm{(well defined for $z$ small enough),}$$ where $U(\cdot)$ is  the analytic function on a neighborhood of zero $U(z)=  \f{-\la-1+\lf[(\la+1)^2+4\la z\ri]^{1/2}}{2\la} $ ($:=z$ if $\la=0$). 

Moreover, by \eqref{17.01.08.2}, we know that $\mu\mapsto C_\mu$ linearizes the rectangular convolution with ratio $\la$: for all $\mu,\nu$ symmetric \pro measures with compact support, \begin{equation}\label{21.1.08.2}C_{\mu\arc\nu}(z)=C_\mu(z)+C_\nu(z)\quad\textrm{($z$ small enough).}\end{equation}

\subsection{Rectangular $R$-transform of unbounded measures: definition}
The  {\textit{Cauchy transform}} and the {\textit{rectangular Cauchy transform with ratio $\la$}}  of symmetric \pro measures with unbounded support are defined with the same formulas as in the case of a compactly supported  measures (see the previous paragraph). 
\begin{propo}\label{Delphine=princessecornichon} Let $A$ be a set of symmetric \pro measures on the real line. Then the following assertions are equivalent

\begin{itemize}
\item[(i)] A is tight,
\item[(ii)] for every $0<\theta < \pi$,  $\underset{\substack{z\to 0\\ \lf|\arg z -\pi\ri|<\theta}}{\lim}\ff{z}H_\mu(z)=1$ uniformly in $\mu\in A$,
\item[(iii)] $\underset{\substack{x\to 0\\ x\in (-\infty,0)}}{\lim}\ff{x} H_\mu(x) =1$  uniformly in $\mu\in A$.
\end{itemize}
\end{propo}

\begin{pr}
$(i)\Rightarrow (ii)$ follows from the well known fact (\cite{defconv}) that if $A$ is tight, $$\forall 0<\theta < \pi, 
\lim_{\substack{z\to 0\\ \lf|\arg z -\f{\pi}{2}\ri|<\f{\theta}{2}}}\ff{z}G_\mu(\ff{z})=1 \quad \textrm{uniformly in $\mu\in A$.}$$  $(ii)\Rightarrow (iii)$ is clear. Suppose $(iii)$. Consider $\eps >0$. 
Take $\eta >0$ \st $\forall u\in [0,1], |\la u^2+(1-\la)u-1|<\eta\Rightarrow |u-1|<\eps/2$. Take $x\in (-\infty,0)$ \st for all $\mu\in A$, $|\ff{x}H_\mu (x)-1|<\eta$. Put $\sqrt{x}=iy$. Then for all $\mu\in A$,  $$ |\ff{iy}G_\mu(\ff{iy})-1|<\eps/2,\textrm{ i.e. }\int_t\f{t^2y^2}{t^2y^2+1}\ud \mu(t)<\eps/2.$$ Let $M$ be \st $\forall t> M, \f{t^2y^2}{t^2y^2+1}\geq \ff{2} $. Then   for all $ \mu\in A$, $\mu([-M,M]^c)\leq 2\int _t\f{t^2{y}^2}{t^2{y}^2+1}\ud \mu(t)<\eps.$
\end{pr}
%Let us denote by $\sqrt{\quad}$ the analytic version of the square root on the complement of the real nonnegative half line in the complex plane \st $\sqrt{-1}=i$.  
%Define, for $0<\alpha < \pi$, $\Gamma_\alpha$ to be the angular sector  of complex numbers $z$ %with negative real part  \st $|\arg z-\pi|< \alpha $. 

Define, for $\alpha \in (0, \pi)$, $\beta > 0$, $\Delta_{\alpha,\beta}$ to be the set of complex numbers $z$ \st $|\arg z-\pi|< \alpha $ and $|z|<\beta$.

Let $\mc{H}$ be the set of functions $f$ which are analytic in a domain $\D_f$ \st for all  $\alpha \in (0, \pi)$, there exists $\beta$ positive \st $$\Delta_{\alpha,\beta}\subset\D_f.$$ A family $(f_a)_{a\in A}$ of functions of $\mc{H}$ is said to be {\it uniform} if for all  $\alpha \in (0, \pi)$,  there exists $\beta$ positive \st $$\forall a\in A,\quad \Delta_{\alpha,\beta}\subset\D_{f_a}.$$

The following \teo has already been  used, in other forms, to define Voiculescu's $R$- and $S$-transforms (see paragraph 5 of \cite{defconv}). Its proof relies on Rouch\'e theorem.

\begin{Th}\label{passiontroisiemeage}
Let $(H_a)_{a\in A}$ be a uniform family of functions of $\mc{H}$ \st for every $\alpha\in (0, \pi)$,  $$\underset{\substack{z\to 0\\ |\arg z-\pi|< \alpha }}{\lim}\f{H_a(z)}{z}=1 \textrm{ uniformly in $a\in A$.}$$

Then there exists a uniform family $(F_a)_{a\in A}$ of functions  of $\mc{H}$ \st for every $\alpha\in (0, \pi)$,  $$\underset{\substack{z\to0\\ |\arg z-\pi|< \alpha }}{\lim}\f{F_a(z)}{z}=1\textrm{ uniformly in $a\in A$,}$$ and there exists $\beta$ positive \st $$\forall a\in A,\quad H_a\circ F_a=F_a\circ H_a=I_d\textrm{ on $\Delta_{\alpha,\beta}$.}$$ 

Moreover, the family $(F_a)_{a\in A}$ is unique in the following sense: if a family $(\tilde{F}_a)_{a\in A}$  of functions  of $\mc{H}$ satisfies the same conditions, then for all  $\alpha \in (0, \pi)$, there exists $\beta$ positive \st $$\forall a\in A,\quad F_a=\tilde{F}_a \textrm{ on }\Delta_{\alpha,\beta}.$$
\end{Th}

For every symmetric \pro measure on the real line $\mu$, let us define the {\textit{rectangular $R$-transform $R_\mu$ with ratio $\la$}}  of $\mu$ by \[C_\mu(z)=U\lf( \f{z}{H_\mu^{-1}(z)}-1\ri),\] where $H_\mu^{-1}$ is defined by the previous \teo  and the function $U$ is the one  defined at the end of section \ref{supremepatience}. 

One can summarize the different steps of the construction of the rectangular $R$-transform with ratio $\la$ in the following chain$$
\begin{array}{l}\ds\underset{\substack{{\textrm{sym. prob.}}\\ {\textrm{measure}}}}{\mu}\,\,\longrightarrow \,\,\underset{\substack{{\textrm{Cauchy}}\\ {\textrm{transf.}}}}{G_\mu}\,\,\longrightarrow\,\, H_\mu(z)=\la G_\mu\lf(\ff{\sqrt{z}}\ri)^2+(1-\la)\sqrt{z}G_\mu\lf(\ff{\sqrt{z}}\ri)\,\,\longrightarrow\,\,\\ \ds  \underset{\textrm{rect. $R$-transf. with ratio $\la$}}{C_\mu(z)=U\lf( \f{z}{H_\mu^{-1}(z)}-1\ri).}\end{array}$$

\subsection{The special cases $\la=0$ and $\la=1$}\label{18.09.06.11}
Note that  the  rectangular $R$-transform with ratio $1$ (resp. $0$), for a symmetric distribution $\mu$,  is linked to  the \trv $\vfi_\mu$ of $\mu$  by the relation $C_\mu(z)=\sqrt{z}\vfi_\mu(1/\sqrt{z})$ (resp. $C_\mu(z)=z\vfi_\rho(1/z)$, where $\rho$ is the push-forward of $\mu$ by the function $t\to t^2$) (see paragraph 5 of \cite{defconv} for the construction of the Voiculescu transform).

\subsection{Remark about the characterization of $R$-transforms} 
The definition of the rectangular $R$-transform of a \pro measure $\mu$ presents many analogies with the definition of its Voiculescu transform $\vfi_\mu$  (see \cite{defconv}). So it seems natural to state, as the authors of \cite{defconv} did for the Voiculescu transform (Proposition 5.6), a characterization of the functions that are the rectangular $R$-transform of a symmetric \pro measure. The characterization of the Voiculescu transform is based on the fact that for every \pro measure $\mu$, the inverse of $z+\vfi_\mu(z)$ extends to  a Pick function (i.e. an analytic function defined on the upper half plane, whose imaginary part does not take negative values), and on the fact that any Pick function equivalent to $z$ at infinity is of the type  $1/G_\mu$, so its inverse is $z+\vfi_\mu(z)$. But unless $\la=0$ or $1$, Pick functions do not appear in an analogous place in the definition of the rectangular $R$-transform, so we cannot proceed similarly to characterize rectangular $R$-transforms. 

Furthermore, the L\'evy Kinchine formula for $\arc$ (see \cite{fbg05.inf.div}), compared with \teo 5.10 of \cite{defconv}) will state that when $\mu$ is $\arc$-infinitely divisible, there exists a unique symmetric \fid $\nu$ \st $C_\mu(z)=\sqrt{z}\vfi_\nu(1/\sqrt{z})$. So the question of the characterization of rectangular $R$-transforms joins another question: can we extend the correspondence $\mu\leftrightarrow \nu$ to a bijective correspondence from the set of symmetric distributions into itself \st for all $\mu$, $C_\mu(z)=\sqrt{z}\vfi_\nu(1/\sqrt{z})$. The analytic functions $f$ on $\gab$ of the type   $f(z)=\sqrt{z}\vfi_\nu(1/\sqrt{z})$, with $\nu$ symmetric \pro measure, are characterized by:
\begin{itemize}\item[(i)] $ f(\bar{z})=\overline{f(z)}$,
\item[(ii)] $\underset{z\to 0}{\lim} f(z)=0$,
\item[(iii)] for all $n$ and all $z_1,\ldots ,z_n$, the matrix 
$\lf[
\f{
\overline{ \sqrt{z_j} }-\sqrt{z_k}
}
{
 \overline{ \sqrt{z_j} }(1+f(z_k)) - \sqrt{z_k}(\overline{ 1+f(z_j) })  
}
\ri]_{k,j=1}^n$ is positive. \end{itemize} 
But nothing allows us to claim that the rectangular $R$-transform of any symmetric distribution satisfies (iii), and that every function that satisfies (i), (ii), and (iii) is the  rectangular $R$-transform of a symmetric distribution.

\subsection{Properties of the rectangular $R$-transform}
\begin{Th}[Injectivity of the rectangular $R$-transform] If the rectangular $R$-transforms with ratio $\la$ of two symmetric \pro measures coincide on a neighborhood of $0$ in $(-\infty,0)$, then the measures are equal.  
\end{Th}

 \begin{pr} If the rectangular $R$-transforms with ratio $\la$ of two symmetric \pro measures $\mu,\nu$ coincide on a neighborhood of $0$ in $(-\infty,0)$, then by uniqueness of analytic continuation, they coincide on a $\Delta_{\alpha,\beta}$, and so $H_\mu=H_\nu$ on this set, and $zH_\mu(z)=zH_\nu(z)$ on this set. So there exists $M>0$ \st for all $y>M$, $$\la ( iy2G_\mu(iy))^2+(1-\la )iyG_\mu(iy)=\la ( iy2G_\nu(iy))^2+(1-\la )iyG_\nu(iy).$$
But if $\rho$ is a symmetric \pro measure, $iyG_\mu(iy)=\int_{t\in \R}\f{y^2}{y^2+t^2}\ud \rho (t)\in (0,1]$. So, by injectivity of $u\mapsto \la u^2+(1-\la )u$ on $(0,1]$, $G_\mu(iy)=G_\nu (iy)$ for $y>1$, and then, by analycity of the Cauchy transform, $G_\mu=G_\nu$, and by injectivity of the Cauchy transform (see \cite{akhi} or section 3.1 of \cite{hiai}) $\mu=\nu$.
\end{pr}

The following remark gives a practical way to derive  any symmetric \pro measure on the real line $\mu$ from $C_\mu$.

\begin{rmq}[How to compute $\mu$ when we know $C_\mu$ ?]\label{How?} First, we have $z/H_\mu^{-1}(z)=T(C_\mu(z))$, for $z\in \C\backslash \R^+$ small enough. From this, we can compute $H_\mu(z)$ for $z\in \C\backslash \R^+$ small enough. Then we can use the equation, for  $z\in \C\backslash \R^+$, $$\ff{z}H_\mu(z)=\la\lf(\ff{\sqrt{z}}G_\mu(\ff{\sqrt{z}})\ri)^2+(1-\la)\ff{\sqrt{z}}G_\mu(\ff{\sqrt{z}}).$$Moreover, when $z\in \C\backslash \R^+$ is small enough, $1/\sqrt{z}$ is large and in $\C^-$, so $\ff{\sqrt{z}}G_\mu\lf(\ff{\sqrt{z}}\ri)$ is closed to $1$. $\ff{z}H_\mu(z)$ is also closed to $1$, and for $h,g$ complex numbers closed to $1$, $$h=\la g^2+(1-\la)g\Leftrightarrow g=V(h),\textrm{ with }V(z)=\f{\la-1+((\la-1)^2+4\la z)^\ff{2}}{2\la}=U(z-1)+1.$$ 
So one has, for $z\in  \C\backslash \R^+$ small enough,  $\ff{\sqrt{z}}G_\mu(\ff{\sqrt{z}})=V(\f{H_\mu(z)}{z}),$ which allows to recover $\mu$ %Note that, for $y>0$, $\sqrt{-1/y^2}=i/y$ and since $\mu$ is symmetric,  $$\ff{\sqrt{-1/y^2}}G_\mu\lf(\ff{\sqrt{-1/y^2}}\ri)=\int\f{y^2}{y^2+t^2}\ud \mu(t).$$
\end{rmq}

The following lemma is an easy consequence of proposition \ref{Delphine=princessecornichon} and \teo \ref{passiontroisiemeage}.
\begin{lem}[Tightness and rectangular $R$-transform]\label{18.11.03.2} Let $A$ be a set of symmetric \pro measures. Then we have equivalence between 

\begin{itemize}
\item[(i)] $A$ is tight,
\item[(ii)] for any $0<\alpha <\pi$,  $\underset{\substack{z\to 0\\ \lf|\arg z -\pi\ri|<\alpha}}{\lim}C_\mu(z)=0$ uniformly in $\mu\in A$,
\item[(iii)] $\underset{\substack{x\to 0\\ x\in (-\infty,0)}}{\lim}R_\mu(x)=0$ uniformly in $\mu\in A$.
\end{itemize} 
\end{lem}

\begin{Th}[Paul L\'evy's \teo for rectangular $R$-transform]\label{tesprof?}  
Let $(\mu_n)$ be a sequence of symmetric \pro measures. Then we have equivalence between:
\begin{itemize}
\item[(i)]  $(\mu_n)$ converges weakly to a symmetric \pro measure;
\item[(ii)] there exists $\alpha,\beta$ \st \begin{itemize}\item[(a)]   $\underset{\substack{z\to 0\\ |\arg z-\pi|<\alpha}}{\lim}C_{\mu_n}(z)=0$ uniformly in $n$, \item[(b)] the sequence $(C_{\mu_n})$ converges uniformly on every compact set of $\Delta_{\alpha,\beta}$;\end{itemize} 
\item[(iii)] \begin{itemize}\item[(a)] $\underset{\substack{x\to 0\\ x\in (-\infty,0)}}{\lim}C_{\mu_n}(x)=0$ uniformly in $n$,  \item[(b)] there exists $\beta> 0$ \st the sequence $(C_{\mu_n})$ converges pointwise  on $[-\beta, 0)$. \end{itemize} 
\end{itemize} Moreover, in this case, denoting by $\mu$ the weak limit of $(\mu_n)$, for every $\alpha$, there exists $\beta$ \st the sequence $(C_{\mu_n})$ converges uniformly to $C_\mu$ on every compact set of $\Delta_{\alpha,\beta}$. 
\end{Th}

\begin{pr} $(i)\Rightarrow (ii)$: suppose that $(\mu_n)$ converges weakly to $\mu$.  Then by the previous lemma, we have $(a)$ of $(ii)$. So there exists $\beta >0$ \st on $\Delta_{\alpha,\beta}$, for all $n$,  $\lf|C_{\mu_n}\ri|\leq 1$. So, by Montel's \teov it suffices to show that the only possible limit, for uniform convergence on every compact, of any subsequence of $(C_{\mu_n})$ is $C_\mu$.  Let $C$ be an analytic function  on $\Delta_{\alpha,\beta}$ \st a subsequence $\lf(C_{\mu_{k_n}}\ri)$ of $(C_{\mu_n})$ converges uniformly to $C$ on every compact of  $\Delta_{\alpha,\beta}$. %By $(a)$ of $(ii)$, 
We know    
(see \cite{akhi} or section 3.1 of \cite{hiai}) that the sequence  $\lf(G_{\mu_n}\ri)$ converges uniformly on every  compact of the upper half plane to $G_\mu$. So the sequence  $\lf(H_{\mu_n}\ri)$ converges uniformly on every  compact of the complement of $[0,+\infty)$ to $H_\mu$. 
Since $\underset{z\to 0}{\lim}C(z)=0$, to prove $C=C_\mu$, it suffices to prove that \begin{eqnarray*}\lf(\la C+1\ri)\lf(C+1\ri)&=&\lf(\la C_\mu +1\ri)\lf(C_\mu+1\ri).\\
{\textrm{But }}\lf(\la C_\mu(z) +1\ri)\lf(C_\mu(z)+1\ri)&=&\f{z}{H_\mu^{-1}(z)}\end{eqnarray*} 
So it suffices to prove that $$\f{z}{\lf(\la C(z) +1\ri)\lf(C(z)+1\ri)}=H_\mu^{-1}(z),$$ so by  \teo \ref{passiontroisiemeage},  it suffices to prove that $$ H_\mu\lf(\f{z}{\lf(\la C(z) +1\ri)\lf(C(z)+1\ri)}\ri)=z.$$
We have \begin{eqnarray*}&\lf|H_\mu\lf[\f{z}{\lf(\la C(z) +1\ri)\lf(C(z)+1\ri)}\ri]-z\ri| &\\
=&\!\!\!\!\!\!\!\!\lf|H_\mu\lf[\f{z}{\lf(\la C(z) +1\ri)\lf(C(z)+1\ri)}\ri]-H_{\mu_{k_n}}\lf[\f{z}{\lf(\la C_{\mu_{k_n}}(z) +1\ri)\lf(C_{\mu_{k_n}}(z)+1\ri)}\ri]\ri| &\\
\leq &\!\!\!\!\!\!\!\!\underbrace{\lf|H_\mu\lf[\f{z}{\lf(\la C(z) +1\ri)\lf(C(z)+1\ri)}\ri]-H_\mu\lf[\f{z}{\lf(\la C_{\mu_{k_n}}(z) +1\ri)\lf(C_{\mu_{k_n}}(z)+1\ri)}\ri]\ri|}_{(1)} &\\ 
&\!\!\!\!\!\!\!\!+\underbrace{\lf|H_\mu\lf[\f{z}{\lf(\la C_{\mu_{k_n}}(z) +1\ri)\lf(C_{\mu_{k_n}}(z)+1\ri)}\ri]-H_{\mu_{k_n}}\lf[\f{z}{\lf(\la C_{\mu_{k_n}}(z) +1\ri)\lf(C_{\mu_{k_n}}(z)+1\ri)}\ri]\ri|}_{(2)} &
\end{eqnarray*} 
By continuity of $H_\mu$, $(1)$ tends to zero when $n$ tends to infinity, and, since $H_{\mu_{k_n}}$ converges uniformly  to $H_\mu$ on every compact, $(2)$ tends to zero. So $H_\mu\lf(\f{z}{\lf(\la C(z) +1\ri)\lf(C(z)+1\ri)}\ri)=z$.

$(ii)\Rightarrow (iii)$ is clear.  

$(iii)\Rightarrow (i)$: suppose $(iii)$. Then by $(i)\Rightarrow (ii)$, every limit of a subsequence of $(\mu_n)$ has a rectangular $R$-transform equal to the pointwise limit of $(C_{\mu_n})$ on $[\beta,0)$. By uniqueness of analytic continuation, all the limits of subsequences of $(\mu_n)$ have the same rectangular $R$-transform, so, by injectivity of this transform (previous theorem), there cannot be more than one  limit of subsequence of $(\mu_n)$. As by $(a)$ of $(iii)$,  the  set $\{\mu_n\ste n\in\n\}$ is tight,  $(\mu_n)$ converges weakly to a symmetric \pro measure.
\end{pr}

\subsection{Rectangular convolution of measures with unbounded support}\label{perez.la.praline} 

\begin{Th}The binary  operation $\arc $ defined on the set of compactly supported symmetric \pro measures in section \ref{ma.raison.vient.tjrs.me.dire}, extends in a unique way to a commutative, associative, and continuous (with respect to the weak convergence) binary operation on the set of  symmetric probability measures on the real line. This operation is called the rectangular free convolution with ratio $\la$. Moreover, for all $\mu,\nu$ symmetric \pro measures one has $C_{\mu\arc\nu}=C_\mu + C_\nu$. \end{Th}

\begin{pr}Let $\mu,\nu$ be symmetric \pro measures. If $(\mu_n)$ (resp. $(\nu_n)$) is a sequence of compactly supported symmetric \pro measures that converges weakly to $\mu$ (resp. $\nu$), then by \teo \ref{tesprof?} and (\ref{21.1.08.2}), the sequence $\lf(\mu_n\arc\nu_n\ri)$ converges to a measure whose rectangular $R$-transform is $C_\mu+C_\nu$ (thus by injectivity of the rectangular $R$-transform, this measure does not depend on the choice of the sequences $(\mu_n)$ and $(\nu_n)$ and is equal to $\mu\arc\nu$ when $\mu$ and $\nu$ are compactly supported).  Note that  $C_{\mu\arc\nu}=C_\mu + C_\nu$ stays true for all symmetric \pro measures $\mu$ and $\nu$. Moreover, this equation shows that $\arc$ is  a commutative, associative, and continuous binary operation on the set of symmetric \pro measures.  \end{pr}

Now we can extend the proposition \ref{15.7.4.1} to the case where $\mu_1,\mu_2$ are not compactly supported: 
\begin{Th}\label{15.7.4.2} Consider sequences $q_1(n),q_2(n)$ \st $$%q_1(n)+q_2(n)=n,\quad 
q_1(n)\leq q_2(n), \quad q_1(n)\ninf +\infty, \quad q_1(n)/\!q_2(n)\ninf \la\in [0,1].$$ Let, for all $n$, $A_n,B_n$ be independent $q_1(n)\tii q_2(n)$ random matrices, one of them being  bi-unitarily invariant,   \st the symmetrizations of the singular  laws of $A_n,B_n$ converge in \pro to \pro measures $\mu,\nu$. Then  the symmetrization of the singular law of $A_n+B_n$ converges in \pro to $\mu\arc\nu$.
\end{Th}

To prove the \teov we will need the following lemma, for which we need to introduce an extended functional calculus. When $F$ is a real  Borel function on the real line, for any  $q_1(n)\tii q_2(n)$ matrix $C$, we  define $F(C)$ to be the matrix $$U[\delta_i^jF(h_i)]_{\substack{1\leq i\leq q_1(n)\\ 1\leq j\leq q_2(n)}}V,$$ where $h_1,\ldots,h_{q_1(n)}\geq 0$, and $U,V$ are respectively $q_1(n)\tii q_1(n)$, $q_2(n)\tii q_2(n)$ unitary matrices \st $$C=U[\delta_i^jh_i]_{\substack{1\leq i\leq q_1(n)\\ 1\leq j\leq q_2(n)}}V.$$ Note that $(F+G)(C)=F(C)+G(C)$ and $F(C)G(C)^*=0$ when $FG=0$. For any \pro measure $\sigma$, we will denote by $F(\sigma)$ (for example $\sigma^2$, $\sigma^{1/\! 2}$, $|\sigma|$\ldots) the push-forward of $\sigma$ by $F$. We denote by $\mu_H$ the spectral measure of an hermitian matrix $H$. Recall that we denote by $\tilde{\sigma}$ the symmetrization of any \pro measure $\sigma$. For example, the symmetrization of the singular law of a rectangular matrix $A$ will be denoted by $\tilde{\mu}_{|A|}$. 

\begin{lem}\label{12.01.03.04}
Let $\lf(C_n\ri)$ be a sequence of $q_1(n)\tii q_2(n)$  random matrices \st $\tilde{\mu}_{|C_n|}$ converges in \pro to a \pro measure $\sigma$. Then, for any odd real function $F$ continuous at $\sigma$-almost  every point of $\R$ and \st $F(\R^+)\subset \R^+$, $\tilde{\mu}_{|F(C_n)|}$ converges in \pro to $F(\sigma)$.   \end{lem}
\begin{pr} $F(\sigma)$ is symmetric and the symmetrization is a continuous operation, so it suffices to prove that $\mu_{|F(C_n)|}$ converges in \pro to $|F(\sigma)|=F(|\sigma|)$. But $\mu_{|F(C_n)|}=F(\mu_{|C_n|})$, $\mu_{|C_n|}$ converges  in  \pro to $|\sigma|$, so the conclusion follows from the fact that the function $\sigma\mapsto F(\sigma)$ on the set of \pro measures on the real line is weakly continuous at $|\sigma|$ because  $F$ is continuous at $|\sigma|$-almost every point of the real line (\cite{billingsley}).\end{pr}

\begin{prth}First of all,  for $U_n$, $V_n$  independent  $q_1(n)\tii q_1(n)$,  $q_2(n)\tii q_2(n)$ uniform random unitary matrices which are independent from $\{A_n,B_n\}$,  the random matrices $U_nA_nV_n, U_nB_nV_n$ are independent (it is a property of the Haar measure on compact groups), have the same singular values as respectively $A_n,B_n$, both of them are   bi-unitarily invariant  and the singular values of $U_nA_nV_n+U_nB_nV_n=U_n(A_n+B_n)V_n$ are the same ones as the ones of $A_n+B_n$. Hence one can suppose that both $A_n,B_n$ are    bi-unitarily invariant.

The function that maps a \pro measure $\sigma$ on  $\R^+$ to the symmetrization $\widetilde{\sigma^{1/\! 2} }$ of its push-forward by $x\to x^{1/\! 2}$ is continuous. So it suffices to prove that the push-forward of the \spl of $|A_n+B_n|$ by $x\to x^2$, i.e. the \spl of $M_n:=(A_n+B_n)(A_n+B_n)^*$ converges in \pro to $(\mu\arc\nu)^2$.  
We can define a distance  on the set of \pro  measures on the real line with the Cauchy transform by $(\sigma_1,\sigma_2)\mapsto \sup\{\lf|G_{\sigma_1} (z)-G_{\sigma_2} (z)\ri|\ste \Im z\geq 1\}.$ This distance defines the topology of weak convergence (\cite{akhi}, \cite{pasturlejay}). 
The Cauchy transform of the spectral distribution of an hermitian matrix $M$ is the normalized trace  of its resolvant $\mf{R}_z(M)=(z-M)^{-1}$. 
So the spectral distribution of a sequence $\lf(X_n\ri)$  of $q_1(n)\tii q_1(n)$  hermitian random matrices converges in \pro to a \pro measure $\sigma$ on the real line \ssi for each $\eps>0$, the \pro of the event $$\underset{\Im z\geq 1}{\sup}\lf|\tr \mf{R}_z\lf(X_n \ri)-G_{\sigma}(z)\ri|>\eps$$ 
tends to zero as $n$ tends to infinity.
Choose  $\eps >0$. We will show that \begin{equation}\label{21.01.08.33}\Pro\{\underset{\Im z\geq 1}{\sup}\arrowvert\tr(\mathfrak{R}_z(M_n))-G_{(\mu\arc\nu)^2}(z)\arrowvert\;>\;\varepsilon\}\stackrel{n\to\infty}{\longrightarrow} 0,\end{equation}where $\Pro$ 
designs the \pro measure of the \pro space where the random matrices are defined. 
Let us define, for $t>0$,  the function $F_t$ on the real line by $$ F_t : x\mapsto \begin{cases}x&{\textrm{ if $|x|\leq t$,}}\\ 0 &{\textrm{ else.}}\end{cases}$$For every \pro measure $\sigma$, $F_t(\sigma)$ converges weakly to $\sigma$ when $t$ tends to infinity, 
so, by continuity of $\arc$, there exists $t\in (0,+\infty)$ \st  $t$ and $-t$ are not atoms of the measures $\mu$ and $\nu$,  and \st \begin{eqnarray}
\underset{\Im z\geq 1}{\sup} | G_{(F_t(\mu)\arc F_t(\nu))^2}(z)-G_{(\mu\arc\nu)^2}(z)|&<&\frac{\varepsilon}{3},\label{12.01.03.02}\\
\mu(\R\bck [-t-1,t+1])+\nu(\R\bck [-t-1,t+1])&<&\frac{\varepsilon}{18}.\label{12.01.03.03}
\end{eqnarray}We will now use the notations $F_t=F$, $G(x)=x-F(x)$ and $M_{n,t}=(F(A_n)+F(B_n))(F(A_n)+F(B_n))^*$.
By triangular inequality and (\ref{12.01.03.02}), we have  \begin{eqnarray}
&\Pro\{\underset{\Im z\geq 1}{\sup}\arrowvert\tr(\mathfrak{R}_z(M_n))-G_{(\mu\arc\nu)^2}(z)\arrowvert\;>\;\varepsilon\}&\nonumber\\
\leq&\Pro\{\underset{\Im z\geq 1}{\sup}\arrowvert\tr[\mathfrak{R}_z(M_n)-\mathfrak{R}_z(M_{n,t})]\arrowvert\;>\;\frac{\varepsilon}{3}\}&\nonumber\\
&+\;\;\;\Pro\{\underset{\Im z\geq 1}{\sup}\arrowvert\tr\mathfrak{R}_z(M_{n,t})-G_{(F(\mu)\arc F(\nu))^2}(z)\arrowvert\;>\;\frac{\varepsilon}{3}\}&\label{12.01.03.05}
\end{eqnarray} 
To treat the first term of the right hand side of (\ref{12.01.03.05}), recall that for any $q_1(n)\tii q_1(n)$ matrix $T$, $|\tr T|\leq \ff{q_1(n)}||T|| \rg T$. The operator norm $||\mathfrak{R}_z(M_n)-\mathfrak{R}_z(M_{n,t})||$  is not less than $2$, because $\Im z\geq 1$. Moreover, due to the equation   $$\mathfrak{R}_z(M_n)-\mathfrak{R}_z(M_{n,t})=-\mathfrak{R}_z(M_n)(M_n-M_{n,t})\mathfrak{R}_z(M_{n,t}),$$ 
the rank of $\mathfrak{R}_z(M_n)-\mathfrak{R}_z(M_{n,t})$  is not more than the one of $M_n-M_{n,t}$. One has (omitting the indices $n$ in $A_n,B_n$) $$M_n-M_{n,t}=G(A)G(A)^*+ G(B)G(B)^*+G(A)B^*+F(A)G(B)^*+G(B)A^*+F(B)G(A)^*,$$ so its rank is not more that $3\rg G(A_n)+3\rg G(B_n)$.
Hence  
\begin{eqnarray*}
&\Pro\{
\underset{\Im z\geq 1}{\sup}
\lf\arrowvert\tr\lf[
\mathfrak{R}_z(M_n)-\mathfrak{R}_z(M_{n,t})\ri]\ri\arrowvert
\;
>
\;
\frac{\varepsilon}{3}
\}
&\\
\leq &\Pro\lf\{\f{6}{q_1(n)}(\rg G(A_n)+\rg G(B_n))>\f{\eps}{3} \ri\}
&\\
=&\Pro\lf\{\mu_{A_n}(\R \bck [-t,t])+\mu_{B_n}(\R \bck [-t,t]) >\f{\eps}{18}\ri\},
&\end{eqnarray*}
which tends to zero when $n$ goes to infinity, by (\ref{12.01.03.03}).
On the other side, the second term of right hand side of (\ref{12.01.03.05}) goes to zero when $n$ goes to infinity by definition  of $\arc$ and by the previous lemma. 
So \eqref{21.01.08.33} is proved.
\end{prth}

\subsection{Examples}
In this section, we give examples of 
rectangular $R$-transforms of symmetric probability measures, and examples of computations of rectangular free convolutions. 

The case where $\la=1$ (see section \ref{18.09.06.11}), where the rectangular $R$-transform and the rectangular free convolution are essentially the Voiculescu transform and the free convolution defined by Voiculescu and Bercovici, doesn't present  any new interest (one can find the little number of  examples known in, e.g., \cite{hiai}), so we will suppose that $\la\in [0,1)$. Unfortunately, for $\la\in (0,1)$, few computations can be done, still less than for the well known case $\la=1$. 

\subsubsection{Convolution of symmetric Bernouilli distributions} 

\begin{propo}\label{8.3.07.stone.roses} Suppose $\la>0$. Then
  $\tau:=\f{\delta_{-1}+\delta_1}{2}\arc\f{\delta_{-1}+\delta_1}{2}$ has support \begin{equation}\label{31.07.07.2}[-(2+\kappa)^{1/2}, -(2-\kappa)^{1/2}]\cup [(2-\kappa)^{1/2},(2+\kappa)^{1/2}],\end{equation}
with $\kappa=2(\la(2-\la))^{1/2}\in (0,2)$, and it admits a density with formula \begin{equation}\label{31.07.07.1}\f{\lf[\kappa^2-(x^2-2)^2\ri]^{1/2}}{\pi\la |x|(4-x^2)}\end{equation} on its support.
\end{propo}

%Before proving this proposition, let us give its concrete signification in the following corollary. 
\begin{rmq}By theorem \ref{15.7.4.2}, this proposition means concretely that for  sequences $q_1(n)\leq q_2(n)$ which tend to infinity \st $q_1(n)/q_2(n)$ tends to $\la \in (0,1)$, if one considers independent random matrices $A_n,B_n,U_n,V_n$ \st  for all $n$, $A_n,B_n$ are $q_1(n)\tii q_2(n)$ matrices \st $A_nA_n^*$, $B_nB_n^*$ are the identity matrix and $U_n$, $V_n$ are uniformly distributed unitary random matrices with respective sizes $q_1(n), q_2(n)$,  then the symmetrized singular law of $A_n +U_nB_n V_n$ tends to the probability measure $\tau$ with the density of \eqref{31.07.07.1} on the support of \eqref{31.07.07.2}. Equivalently, one can say that the spectral law of $$(A_n +U_nB_n V_n)(A_n +U_nB_n V_n)^*$$ tends to the push forward of $\tau$ by the map $x\mapsto x^2$, i.e. to the probability measure with support $[2-\kappa,2+\kappa]$ and with density $$\f{\lf[\kappa^2-(x-2)^2\ri]^{1/2}}{\pi\la x(4-x)}.$$\end{rmq}

\begin{rmq}Note that the rectangular free convolution of the Dirac mass at zero with any other symmetric \pro measure is the measure itself, and that for any $c>0$, $\f{\delta_{-c}+\delta_c}{2}\arc\f{\delta_{-c}+\delta_c}{2}$ can be deduced from this  proposition using a dilation. Note also that if $\la=0$, then $\tau=\f{\delta_{-\sqrt{2}}+\delta_{\sqrt{2}}}{2}$. \end{rmq}

\begin{pr}Define $\nu=\f{\delta_{-1}+\delta_1}{2}$. We have
$G_\nu(z)= \f{z}{z^2-1},$ hence $G_\nu(1/\sqrt{z})=\f{\sqrt{z}}{1-z}$ and 
$$H_\nu(z)=\f{\la z}{(1-z)^2}+(1-\la)\f{z(1-z)}{(1-z)^2}=\f{(\la-1)z^2+z}{(1-z)^2}.$$ Hence for $x,y\in \C\bck\R^+$ in a neighborhood of zero, the equation $H_\nu(y)=x$ is equivalent to $(\la-1)y^2+y=x-2xy+xy^2$, i.e. to $y^2(x+1-\la)-y(1+2x)+x=0.$ %The discriminent of this second degree polynomial is $$\Delta=(1+2ax)^2-4x(a^2x+(1-\la)a)=1+4ax+4a^2x^2-4a^2x^2-4(1-\la)ax=1+4\la ax.$$
In this case, since $y=H_\nu^{-1}(x)$ must tend to zero as $x$ tends to zero, one has
$$H_\nu^{-1}(x)=\f{1+2x-(1+4\la x)^\ff{2}}{2(x+1-\la)}\textrm{ and }\f{x}{H_\nu^{-1}(x)}-1=\f{(1+4\la x)^\ff{2}-1+2x}{2}.$$
Note that $(\la+1)^2+4\la \lf(\f{x}{H_\nu^{-1}(x)}-1\ri)=\lf(\la+(1+4\la x)^\ff{2}\ri)^2,$ hence $$C_\nu(x)=\f{-\la-1+\lf[  (\la+1)^2+4\la \lf(\f{x}{H_\mu^{-1}(x)}-1\ri)\ri]^\ff{2}}{2\la}=\ff{2\la}\lf[(1+4\la x)^\ff{2}-1\ri].$$
 Thus by \eqref{21.1.08.2}  and  remark \ref{How?},  for $y\in \C\bck\R^+$
in the  neighborhood of zero, if one denotes $(1+4\la y)^\ff{2}$ by $w$, one has  $$H_\tau^{-1}(y)=\f{\la y}{1+4\la y+(\la-1)(1+4\la y)^\ff{2}}=\f{(w^2-1)/4}{w^2+(\la-1)w}.$$Hence for $x\in \C\bck\R^+$
in the  neighborhood of zero,  if one denotes $H_\tau(x)$ by $y$ and $(1+4\la y)^\ff{2}$ by $w$, one has $H_\tau^{-1}(y)=x$, i.e. $w^2-1=4xw^2+4x(\la-1)w,$ i.e., since $w$ tends to $1$ as $x$ tends to zero,  $$w=\f{1}{1-4x}[2x(\la-1)+[1-4x+4(\la-1)^2x^2]^{\ff{2}}].$$ Then an easy computation leads to 
$$(\la-1)^2+\f{4\la y}{x}=\lf(\f{
\la-1+2[1-4x+4(\la-1)^2x^2]^{\ff{2}}}{1-4x}
\ri)^2.$$
Thus by remark \ref{How?}, for $z\in \C\bck\R^+$
in the  neighborhood of zero, $$\ff{\sqrt{z}}G_\tau\lf(\ff{\sqrt{z}}\ri)=V\lf(\f{H_\tau(z)}{z}\ri)=\ff{2\la}\lf(\la-1+\f{
\la-1+2[1-4z+4(\la-1)^2z^2]^{\ff{2}}}{1-4z}\ri),$$ hence $G_\tau\lf(\ff{\sqrt{z}}\ri)=\ff{2\la\ff{\sqrt{z}}}\lf(\la-1+\f{
\la-1+2[1-4/(1/\sqrt{z})^2+4(\la-1)^2/(1/\sqrt{z})^4]^{\ff{2}}}{1-4\ff{\lf(1/\sqrt{z}\ri)^2}}\ri),$ so for $z\in i\R^-$ large enough, $G_\tau(z)=\f{(2z^2-4)(\la-1)+2\lf[(z^2-2)^2-\kappa^2\ri]^{1/2}}{2\la z(z^2-4)}$
and by analytic continuation, for any $z\in \C^-$,  $G_\tau(z)=\f{(4-2z^2)(\la-1)+2s(z)}{2\la z(4-z^2)}$, where $s(z)$ is analytic on $\C^-$ and satisfies $s(z)^2=(z^2-2)^2-\kappa^2$.
  In order to compute $\tau$, we shall use the following well known result (lemma 2.17 of \cite{BA}): %for the two first ones, and lemma 4.7 of \cite{fbg05.inf.div} for the last one).
  for Lebesgue-almost all real number  $x$, the limit, as  $z=u+iv\in \C^-$ tends to $x$ in such a way that $(u-x)/v$ stays bounded, of $\Im G_\tau(z)/\pi$ is the density of the absolutely continuous part of $\tau$ at $x$.
  
%  (b) For all real number $x_0$, $\tau(\{x_0\})=\ds\lim_{z\to x_0}(z-x_0)G_\tau(z)$.
  
  %(c) For $I$ is an open interval, if the restriction of $G_\tau$ to $\C^-$ extends analytically to  an open set containing $\C^-\cup I$, then the restriction of $\tau$ to $I$ admits an analytic density: $x\in I \mapsto \ff{\pi}\Im G_\tau(x).$

  Let us define $$O:=(-(2+\kappa)^{1/2}, -(2-\kappa)^{1/2})\cup ((2-\kappa)^{1/2},(2+\kappa)^{1/2}).$$ By what precedes, for Lebesgue-almost all  $x\in \R$, there is no sequence   $w_n>0$ tending to $0$ \st   $\ds\lim_{n\to\infty}\Im G_\tau(x-iw_n)<0$. Hence since for all $x\in O$, $(x^2-2)^2-\kappa^2$ is negative,  for Lebesgue-almost all  $x\in O$, there is $\eps>0$ \st for all $w\in (0,\eps)$, for $z=x-iw$,$$G_\tau(z)=\begin{cases}
  \f{(4-2z^2)(\la-1)+2\sqrt{(z^2-2)^2-\kappa^2}}{2\la z(4-z^2)}&\textrm{if $x>0$,}\\
  \f{(4-2z^2)(\la-1)-2\sqrt{(z^2-2)^2-\kappa^2}}{2\la z(4-z^2)}&\textrm{if $x<0$.}
  \end{cases}$$
  It implies that for Lebesgue-almost all  $x\in O$, $$\ds\lim_{y\to 0^-} \Im G_\tau(x+iy
 )/\pi =\f{\Im(\sqrt{(x^2-2)^2-\kappa^2})}{\pi\la |x|(4-x^2)}= \f{\lf[\kappa^2-(x^2-2)^2\ri]^{1/2}}{\pi\la |x|(4-x^2)}.$$
 
 Now, to complete the proof, it suffices to prove that $\tau$ is equal the absolutely continuous part of its restriction to $O$, i.e. that $$\int_O \f{\lf[\kappa^2-(x^2-2)^2\ri]^{1/2}}{\pi\la |x|(4-x^2)}\ud x =1.$$
 This can be proved using the changing of variable  $y=\f{x^2-2}{(\kappa^2-(x^2-2)^2)^{1/2}}$ and the well known formula of the Cauchy transform of the Cauchy distribution.
\end{pr}

\subsubsection{The rectangular analogues of the Gaussian distributions} Here, $\la>0$ (the case $\la=0$ can be found in  \cite{fbg05.inf.div}).
It is proved in \cite{fbg05.inf.div} that for all $\sigma>0$, the rectangular analogue of the symmetric Gaussian distribution with variance $\sigma^2$ is the law $N_{\sigma^2}$ with density $$\f{\lf[4\la -(\sigma^2x^2-1-\la)^2\ri]^{1/2}}{2\pi \la \sigma^2 |x| }\chi(\sigma^2x^2),$$where $\chi$ stands for the characteristic function of the interval $[(1-\la^{1/2})^2,(1+\la^{1/2})^2]$. 
 Its support is $[-\sigma(1+\la^{1/2}),-\sigma(1-\la^{1/2})]\cup [\sigma(1-\la^{1/2}),\sigma(1+\la^{1/2})] $, its rectangular $R$-transform is $\sigma^2 z$. 
One has, for all $\sigma, \alpha>0$, $N_{\sigma^2}\arc N_{\alpha^2}=N_{\sigma^2+\alpha^2}$. 

\begin{rmq}[Growth of the support in the related semigroup]\label{growth.support.21.09.06.normal} Note that in the additive semigroup $(\{N_{c}\ste c>0\}, \arc)$ the size of the support of $N_c$ in not linear in $c$ but in $c^\ff{2}$. This kind of phenomenon, which recalls us that we are in a gaussian context, but still surprising since the support (and not the variance) is concerned, had already been observed in free probability.\end{rmq}

\subsubsection{The rectangular analogues of the Cauchy distributions} 

 For all $t>0$, the rectangular analogue of the symmetric Cauchy law with parameter $t$ by the  Bercovici-Pata bijection with ratio $\la$ is $$\nu_t:=1_{\R\backslash [-\f{t(1-\la)}{2}, \f{t(1-\la)}{2}]}(x) \f{t}{ \pi(\la t^2+x^2)}\lf[1-\f{t^2(\la-1)^2}{4x^2}\ri]^\ff{2}\ud x.$$Its support is $\R\backslash (-\f{t(1-\la)}{2}, \f{t(1-\la)}{2})$.
It is proved in \cite{fbg05.inf.div} that the $R$-transform with ratio $\la$ of $\nu_t$ admits an analytic extension to $\C\backslash \R^+$ given by the formula $it\sqrt{z}$. Hence one has, for all $t_1, t_2>0$, $\nu_{t_1}\arc \nu_{t_2}=\nu_{t_1+t_2}$. 
 
\subsubsection{The rectangular analogues of other symmetric stable laws} It can also be proved that for all $\alpha\in (0,2)$, the rectangular analogue of the classical semigroup   of stable laws with index $\alpha$ is the set of laws with rectangular $R$-transforms with ratio $\la$ given by the formula $z\mapsto -C(-z)^{\alpha/2}$, with $C>0$ and where the $\alpha/2$-th power is defined with respect to the canonical determination of the argument on 
 $\C\bck \R^-$. Unfortunately,  though we give their rectangular $R$-transforms with ratio $\la$, except for $\alpha=1$ (see the previous paragraph), we cannot give the densities of these laws. However, just like for free stable laws (see \cite{appenice}), results can be proved about their densities and supports (see \cite{BBG07}).

 \subsubsection{The rectangular analogues of the symmetric Poisson distributions}\label{hetriessospeaking} It is proved, in \cite{fbg05.inf.div}, that for all $c>0$, the sequence  $$\lf(\lf(1-\f{c}{n}\ri)\delta_0+\f{c}{2n}(\delta_{-1}+\delta_1)\ri)^{\arc n}$$has a weak limit, that we shall denote by $P_c$, and that 
the rectangular $R$-transform of   $P_c$  is given by the formula  $\f{cz}{1-z}.$ We cannot compute $P_c$ itself, but 
 we know that for all $c,c'>0$, $P_c\arc P_{c'}=P_{c+c'}$.

 \section{Appendix: a technical result}
In this section, we do not make any particular hypothesis about $\rho_k$'s: they are nonnegative, and at most one of them is zero. So we can suppose that $\rho_2,\ldots, \rho_d$ are positive.

For all integers $n,k$, for each partition $\mc{P}$ of $[k]$, $[n]^\mc{P}$ will design the set of elements $i=(i_1,\ldots,i_k)$ of $[n]^{k}$ such that  $h\stackrel{\mc{P}}{\sim}h'\Leftrightarrow i_h=i_{h'}$. Moreover, $|\mc{P}| $ denotes the number of blocks of $\mc{P}$.

\begin{Def} An $n\tii n$ matrix is said to be an {\rm almost diagonal matrix} if it is of the type $\widetilde{M}$, with $M$   diagonal (either square or rectangular) matrix. In this case, for all $j\in [n]$, $A_j$ will denote the only nonzero term on the $j$-th column of $A$  if it exists, and zero in the other case.
\end{Def}

\begin{lem}\label{01.11.03.1} Let $K$ be a positive integer, and let for each $n\geq 1$, $A(1,n),\ldots ,A(K,n)$ be a family of almost diagonal $n\tii n$ matrices. Let $\mc{H}$ be a subset of $[K]$, and  let $\mc{P}$  be a partition of $[K]$  \st for all $\{h\}$ singleton class of $\mc{P}$, $A(h,n)$ is a diagonal matrix and $\Tr A(h,n)=0$. We suppose too that the family $(A(i,n))_{i,n}$ is uniformly bounded, and that for all $h\in\mc{H}$, for all $i\in \{q_1(n)+1,\ldots, n\}$, $A_{i}(h,n)=0$.
Then$$\sum_{i\in [n]^\mc{P}}\prod_{h=1}^{K}A_{i_h}(h,n)=O(q_1(n)^{|\mc{P}_\mc{H}|-p_\mc{H}/2}{n}^{|\overline{\mc{P}}_\mc{H}|-\overline{p}_\mc{H}/2}),$$
where  $\mc{P}_\mc{H}=\{B\cap\mc{H}\ste B\in \mc{P}\}\bck\{\emptyset\} $ and  $\overline{\mc{P}}_\mc{H}\{B\cap ( K\bck\mc{H})\ste B\in \mc{P}\}\bck\{\emptyset\}$)  and $p_\mc{H}$ (resp. $\overline{p}_\mc{H}$) the number of singletons in $\mc{P}_\mc{H}$ (resp. of $\overline{\mc{P}}_\mc{H}$).  
\end{lem}

\begin{pr}
We will prove this result by induction on $p=p_\mc{H}+\overline{p}_\mc{H}$. If $p=0$, the result is clear. Now, let us suppose the result to be true to the ranks $0,\ldots ,p-1$, with $p\geq 1$. Consider $\mc{P}$ such that $p_\mc{H}+\overline{p}_\mc{H}=p$, and a singleton class $\{h_0\}$ of $\mc{P}$.  
Denote by $\mc{P}'$ the partition $\mc{P}-\{\{ h_0 \}\}$ of $\{1,\ldots, \hat{h}_0, \ldots, K\}$ and let, for each  class $B$ of $\mc{P}$ with $B\neq \{ h_0 \}$,  $\mc{P}(h_0\to B)$ be the partition of $[K]$ obtained from $\mc{P}$ by linking the classes $\{h_0\}$ and $B$. We have  % where the class $\{ h_p \}$ has been linked to the class $B$. 
%We have $|\mc{P}_B|=|\mc{P}|-1$ and the number of singleton classes of $\mc{P}_B$ that correspond to traceless matrices is $p-1$ if $B$ is not in one of these classes in $\mc{P}$, and $p-2$ in the other case.
  %\begin{eqnarray*}&
$$\ds\sum_{i\in [n]^\mc{P}}\prod_{h=1}^K A_{i_h}(h,n)%&\\
=\ds\sum_{(i_1,\ldots,\hat{i}_{h_0},\ldots,i_K)\in [n]^{\mc{P}'}}\lf[\lf(\sum_{i\in[n]-\{i_1,\ldots,\hat{i}_{h_0},\ldots,i_K\}}A_{i}(h_0,n)\ri)%&\\
%&\ds\times
\prod_{\substack{h=1\\ h\neq h_0}}^KA_{i_h}(h,n)\ri].$$%&\end{eqnarray*}
But, since $A(h_0,n)$ is diagonal and has null trace,  $$\ds\sum_{i\in[n]-\{i_1,\ldots,\hat{i}_{h_0},\ldots,i_K\}}A_{i}(h_0,n)=-\sum_{i\in\{i_1,\ldots,\hat{i}_{h_0},\ldots,i_K\}}A_{i}(h_0,n).
$$
%\begin{eqnarray*}&\sum_{\alpha\in[d']\backslash\{\alpha_1,\beta_1,\ldots,\hat{\alpha}_{r_p},\hat{\beta}_{r_p},\ldots,\alpha_L,\beta_L\}}D(r_p,d)_{\alpha,\alpha}&\\
%&=-\sum_{\alpha\in\{\alpha_1,\beta_1,\ldots,\hat{\alpha}_{r_p},\hat{\beta}_{r_p},\ldots,\alpha_L,\beta_L\}}D(r_p,d)_{\alpha,\alpha}.&\end{eqnarray*}
So $$\ds    \sum_{i\in [n]^\mc{P}}\prod_{h=1}^K A_{i_h}(h,n)     
= -\!\!\sum_{\substack{B\in\mc{P}\\ B\neq \{h_0\}}}\sum_{i\in [n]^{\mc{P}(h_0\to B)}}\prod_{h=1}^K A_{i_h}(h,n).$$
For example, suppose  $h_0\in \mc{H}$ (the other case is treated in the same way). With the induction hypothesis, and dividing the sum in $B\in \mc{P}_\mc{H}$ and $B\in \overline{\mc{P}}_\mc{H}$, one has (each $\epsilon_B$ being $1$  or $0$ according to whether $B$ is a singleton or not)  %\begin{eqnarray*}&\ds    \sum_{i\in [d']^\mc{P}}\prod_{h=1}^K A_{i_hi_h}(h,d)     &\\=&\ds -\sum_{\substack{B\in\mc{P}\\ B\neq \{h_p\}}}O(d^{|\mc{P}|-1-(p-\varepsilon_B)/2}.& \end{eqnarray*}
\begin{eqnarray*}
\ds    \sum_{i\in [n]^\mc{P}}
\prod_{h=1}^K A_{i_h}(h,n)     
&= &\ds  \sum_{
\substack{B\in\mc{P}_\mc{H}\\ B\neq \{h_0\}}}
O\lf(
q_1(n)^{
|\mc{P}_\mc{H}|-1-(p_\mc{H}-1-\epsilon_B)/2
}
{n}^{
|\overline{\mc{P}}_\mc{H}|-\overline{p}_\mc{H}/2
}  
\ri)\\ &&\ds
+\!\!\sum_{
\substack{B\in\overline{\mc{P}}_\mc{H}}}
O\lf(q_1(n)^{
|\mc{P}_\mc{H}|-(p_\mc{H}-1)/2}{n}^{
|\overline{\mc{P}}_\mc{H}|-1-(\overline{p}_\mc{H}-\epsilon_B)/2}  
\ri)\\
&= &\ds q_1(n)^{\f{\epsilon_B-1}{2}}\!\!\sum_{
\substack{B\in\mc{P}_\mc{H}\\ B\neq \{h_0\}}}
O\lf(
q_1(n)^{
|\mc{P}_\mc{H}|-p_\mc{H}/2
}
{n}^{
|\overline{\mc{P}}_\mc{H}|-\overline{p}_\mc{H}/2
}  
\ri)\\ &&\ds +q_1(n)^{\ff{2}}{n}^{\f{\epsilon_B-2}{2}}\!\!
\sum_{
\substack{B\in\overline{\mc{P}}_\mc{H}}}
O\lf(q_1(n)^{
|\mc{P}_\mc{H}|-p_\mc{H}/2}{n}^{
|\overline{\mc{P}}_\mc{H}|-\overline{p}_\mc{H}/2}  
\ri). \end{eqnarray*}
Since $q_1(n)=O(n)$, the result is proved. 
\end{pr}

The following proposition, main result of this subsection, is divided in two cases, which give very similar results with very similar proofs.  The proposition involves  special classes of matrices. %An {\it almost diagonal matrices} is an $n\tii n$ matrix whose blocks are all diagonal (either square or rectangular) matrices. 
\begin{Def} For $k\in [d]$, we define $\mathbb{U}_k(n)$ to be the set of $n\tii n$ matrices whose $(k,k)$-th block is a $q_k(n)\tii q_k(n)$ unitary  matrix, and with other blocks zero.\end{Def} 
Note that $(\mathbb{U}_k(n),\times)$ is a compact group, isomorphic to the one of $q_k(n)\tii q_k(n)$  unitary  matrices, and random matrices of $\mathbb{U}_k(n)$ distributed according to the Haar measure will be said to be {\it uniform.}       In the following, we are going to use an integer $N$ and the cyclic order on $[N]$. This means that to put the index $N+1$ on an element is equivalent to put the index $1$. 

\begin{propo}\label{TX1(n).merle.com}Let, for $n\geq 1$, $V(s,k,n)$ ($s\in\n,k\in [d]$), be a family of independent random matrices, \st for all $s,k$, $V(s,k,n)$ is uniform on  $\mathbb{U}_k(n)$. Fix $N$ positive integer and $R>0$. 

\uo Let, for each $n\in\n$, $D(1,n),\ldots,D(N,n)$ be $n\tii n$ constant matrices \st for all $r\in [N]$, $||D(r,n)||\leq R$, and there exists $u(r),v(r)\in [d]$ \st one of the two following conditions is realized: 
\begin{itemize}
\item[(i)] $u(r)=v(r)$ and for all $n$, $D(r,n)=p_{u(r)}(n)$,
\item[(ii)] for all $n$, $D(r,n)$ is an almost diagonal matrix of the type $\widetilde{M}$, where $M$ is a $q_{u(r)}(n)\tii q_{v(r)}(n)$ diagonal matrix, and $\ED(D(r,n))\ninf 0$.   
\end{itemize}
Consider $(s_l,k_l)_{l\in [N]}\in (\n\tii [d])^N$, and suppose that  for all $r\in [N]$ \st (i) above is satisfied, $(s_{r},k_r)\neq (s_{r+1},k_{r+1})$. Consider also $m_1,\ldots,m_N\in \z-\{0\}$ and $\eta>0$. Then  the probability of the event $$\{|| \ED(V(s_1,k_1,n)^{m_1}D(1,n)\cdots D(N-1,n)V(s_N,k_N,n)^{m_N}D(N,n))||\leq \eta\}$$ tends to $1$ as $n$ goes to infinity.

\deuzio  Let, for each $n\in\n$, $D(0,n),\ldots,D(N,n)$ be $n\tii n$ constant matrices with norms also $\leq R$. We suppose that  there exists $v(0), u(N)\in [d]$ \st $D(0,n)$ (resp. $D(N,n)$) is an almost diagonal matrix of the type $\widetilde{M}$, where $M$ is a $q_{1}(n)\tii q_{v(0)}(n)$ (resp. $q_{u(N)}(n)\tii q_{1}(n)$), and that  for all $r\in [N-1]$, there exists $u(r),v(r)\in [d]$ \st one of the two following conditions is realized: 
\begin{itemize}
\item[(i)] $u(r)=v(r)$ and for all $n$, $D(r,n)=p_{u(r)}(n)$,
\item[(ii)] for all $n$, $D(r,n)$ is an almost diagonal matrix of the type $\widetilde{M}$, where $M$ is a $q_{u(r)}(n)\tii q_{v(r)}(n)$ diagonal matrix, and $\ED(D(r,n))\ninf 0$.   
\end{itemize}
Consider $(s_l,k_l)_{l\in [N]}\in (\n\tii [d])^N$, and suppose that  for all $r\in [N-1]$ \st (i) above is satisfied, $(s_{r},k_r)\neq (s_{r+1},k_{r+1})$. Consider also $m_1,\ldots,m_N\in \z-\{0\}$ and $\eta>0$. Then  the probability of the event $$\{|| \ED(D(0,n)V(s_1,k_1,n)^{m_1}D(1,n)\cdots D(N-1,n)V(s_N,k_N,n)^{m_N}D(N,n))||\leq \eta\}$$ tends to $1$ as $n$ goes to infinity.
\end{propo}

\begin{rmq}\label{TX1(n).merle.com'} By linearity of $\ED$, the result stays true in case \uov (resp. in case $2^\circ)$) if  for all $r$ in $[N]$  (resp. in $[N-1]$) \st  (ii) is satisfied, $B(r,n)$ is replaced by a sum of constant almost diagonal matrices whose images by $\ED$ also tend  to zero.\end{rmq}

The proof is inspired from the one of Theorem 4.3.1 p. 147 in \cite{hiai}.

\begin{pr} Both proofs will be made together, we will only have to separate them sometimes. 

{\it Step I. } \udeuzio Let us %prove the result by induction on $N$. If $N=0$, it is clear. So suppose the result to be proved to the rank $N-1$, and let us prove it to the rank $N$.
%Let us consider, for all $r\in [N]$, $u(h),v(h)\in     [d]$ \st we have not $B(h,n)=I_n$ for all $n$, we have: $u(h),v(h)\in     [d]$ \st for all $n$, $B_h(n)$ has only one nonzero block, which is the $(u(h),v(h))$-th. Then 
denote, 
for all $r$ which satisfies (ii) above and  $u(r)=v(r)$,  $$D(r,n) = D'(r,n)+\la_r(n)p_{u(r)}(n), \textrm{ with } \ED(D'(r,n))=0.$$ Then by hypothesis, $\la_h(n)\ninf 0$. Hence, by linearity and boundedness hypothesis, it suffices to prove the result for $D(r,n)$ replaced by $D'(r,n)$. Thus we can from now on suppose that for all $r$ \st (ii) is realized, we have $\ED(D(r,n))=0$ for all $n$.  

 {\it Step II. }  %Let us first prove that we can remove $B_0(n)$ and that it suffices to prove the boundedness of the expectation of \begin{equation}\label{18.1.05.29}\Tr V(s_1,k_1,n)^{m_1}B_1(n)\cdots B_N(n)\Tr B_N(n)^*V(s_N,k_N,n)^{-m_N}B_{N-1}(n)^*\cdots V(s_1,k_1,n)^{-m_1}.  \end{equation} 
\uo By definition of $\ED$, it suffices to prove that for all $k\in [d]$, the normalized trace of the $k$-th diagonal block of the product $$V(s_1,k_1,n)^{m_1}D(1,n)\cdots D(N-1,n)V(s_N,k_N,n)^{m_N}D(N,n)$$converges in \pro to zero, i.e. that \begin{equation}\label{18.1.05.28}\ff{q_k(n)}\Tr p_k(n)V(s_1,k_1,n)^{m_1}D(1,n)\cdots D(N-1,n)V(s_N,k_N,n)^{m_N}D(N,n)p_k(n)\end{equation} converges in \pro to zero. So let us fix $k \in [d]$, and let us prove it. 

If $k_1\neq k$ or $v(N)\neq k$, then the $k$-th block of the matrix is always zero, so the result is clear. So let us suppose that $k_1=v(N)=k$. We can then remove $p_k(n)$ in (\ref{18.1.05.28}).
By Markov inequality, it suffices to prove that (\ref{18.1.05.28}) tends to zero in $L^2$. Since for all matrix $M$, $|\Tr M|^2=\Tr M\Tr M^*$, we only have to prove that the expectation of  \begin{equation}\label{18.1.05.29} \Tr V(s_1,k_1,n)^{m_1}D(1,n)\cdots D(N, n)\Tr D(N,n)^*V(s_N,k_N,n)^{-m_N}D({N-1},n)^*\cdots V(s_1,k_1,n)^{-m_1}\end{equation}  is $o(q_k(n)^2)$.
%{\it Step II. } Let us prove that we can suppose that for all $l\in [N]$, if  we have not $\ED(B_l(n))\neq I_n$ for all $n$, then for all $n$, $\ED(B_l(n))=0$.\\We first prove it for $l<N$, and then for $l=N$.\\
%Let us denote, for all   $l\in [N-1]$, \st we have not $\ED(B_l(n))\neq I_n$ for all $n$,  $B_l(n)=B'_l(n)+\sum_{k=1}^d\la_k(l,n)p_k(n)$, with $\ED(B'_l(n)=0$. Then for all $k,l$, $\la_k(l,n)\ninf 0$. So by linearity  and by boundedness hypothesis, it suffices to prove the result for $B_l(n)$'s replaced by $B'_l(n)$'s. Thus we can from now on suppose that for all such $l$, for all $n$, $\ED(B_l(n))=0$. \\
%In the same way, let us denote $B_N(n)=B'_N(n)+\sum_{k=1}^d\la_k(n)p_k(n)$. 

\deuzio We have to prove that \begin{equation*}%\label{21.1.05.28}
\ff{q_1(n)}\Tr D(0,n)V(s_1,k_1,n)^{m_1}D(1,n)\cdots D(N-1,n)V(s_N,k_N,n)^{m_N}D(N,n)\end{equation*} converges in \pro to zero. Since $\Tr XY=\Tr YX$, it suffices to prove that \begin{equation}\label{21.1.05.28}\ff{q_1(n)}\Tr V(s_1,k_1,n)^{m_1}D(1,n)\cdots D(N-1,n)V(s_N,k_N,n)^{m_N}D(N,n)D(0,n)\end{equation} converges in \pro to zero. From now on, we will denote $D(N,n)D(0,n)$ by $D(N,n)$, which is now an $n\tii n$ matrix of the type $\widetilde{M}$ with $M$ $q_{u(N)}(n)\tii q_{v(0)}(n)$, but which has not more than $q_1(n)$ nonzero entries. For the same reason as above, we only have to prove that  the expectation of \begin{equation}\label{21.1.05.29}\Tr V(s_1,k_1,n)^{m_1}D(1,n)\cdots D(N,n)\Tr D(N,n)^*V(s_N,k_N,n)^{-m_N}D(N-1,n)^*\cdots V(s_1,k_1,n)^{-m_1}\end{equation}  is  $o(q_1(n)^2)$.

 {\it Step III. } \udeuzio Let us expand both traces in the previous product, and re index the sum on partitions.

Let $M=|m_1|+\cdots +|m_N|$, let $M(0)=0$ and, for $r\in \{1,\ldots,N\}$,$$M(r)=|m_1|+\cdots +|m_r|,$$  let $M(N+1)=2M$, and let,
for $r\in\{N+2, \ldots, 2N\}$, $$M(r)=M+|m_N|+|m_{N-1}|+\cdots +|m_{2N+2-r}|.$$ Let also, for $r\in\{N+1,\ldots,2N\}$,  $D(r,n)=D(2N+1-r,n)^*$.
We write, for $\kappa\in \n\tii [d]$, \begin{eqnarray*}v_{i,j}(\kappa,1,n)&=&V(\kappa,n)_{i,j},\\   v_{i,j}(\kappa,-1,n)&=&\bar{V}(\kappa,n)_{j,i}.\end{eqnarray*} With those notations, there exists two functions $\kappa$ and $\eps$ \st the expectation of (\ref{18.1.05.29}) is 

$$\ds\sum_{i\in[n]^{2M}}\;\;\prod_{r=1}^{2N}D(r,n)_{ j_{M(r)}  }\E\lf[ \prod_{h=1}^{2M}v_{i_h,j_h}(\kappa(h),\eps(h),n) \ri],$$ where $\E$ denotes expectation and for all $i=(i_1,\ldots, i_{2M})\in [n]^{2M}$, $$%j_0=i_1,\,
j_1=i_2,\,j_2=i_3,\ldots, j_M=i_1,\;\; j_{M+1}=i_{M+2},\,j_{M+2}=i_{M+3},\,\ldots, j_{2M-1}=i_{2M},\, j_{2M}=i_{M+1}.$$

Let us introduce the partition $\mc{Q}$ of $[2M]$ defined by  $h\stackrel{\mc{Q}}{\sim}h'$ \ssi $i_h=i_{h'}$. %and $(h,h')\in H^2\cup ([2k]\backslash H)^2$, 
Then we can rewrite the preceding sum 
$$\ds \sum_{\substack{\mc{Q}\textrm{ partition}\\ \textrm{of $[2M]$}}}\;\;\sum_{\substack{i\in [n]^{\mc{Q}}}}\;\;
\prod_{r=1}^{2N}D(r,n)_{ j_{M(r)}  }\E\lf[ \prod_{h=1}^{2M}v_{i_h,j_h}(\kappa(h),\eps(h),n) \ri].
$$
Thus it suffices to prove that for all partition $\mc{Q}$ of $[2M]$, the sum  
\begin{equation}\label{emouvantgrosamericain'}\ds  \sum_{\substack{i\in [n]^{\mc{Q}}}}\;\;
\prod_{r=1}^{2N}D(r,n)_{ j_{M(r)}  }\E\lf[ \prod_{h=1}^{2M}v_{i_h,j_h}(\kappa(h),\eps(h),n) \ri]\end{equation}
 is  $o(q_k(n)^2)$ in case \uov and $ o(q_1(n)^2)$ in case $2^\circ$),  as $n$ goes to infinity.

So we fix a partition $\mc{Q}$ of $[2M]$. 

 {\it Step IV. } \udeuzio Let us denote, for all $h\in [2M]$, $\kappa(h)=(s(h), k(h))$. Note that, for $i\in [n]^\mc{Q}$, for \begin{equation}\label{21.1.05.1}\E\lf[ \prod_{h=1}^{2M}v_{i_h,j_h}(\kappa(h),\eps(h),n) \ri]
\end{equation}
to be nonzero, $i$ has to satisfy  \begin{equation}\label{19.1.05.17.35}\forall h\in [2M], i_h,j_h\in \{q_1(n)+\cdots +q_{k(h)-1}(n)+1,\ldots,  q_1(n)+\cdots +q_{k(h)-1}(n)+q_{k(h)}(n)\}.\end{equation} Indeed, in the other case, one of the factors in the product is zero, by definition of $M\mapsto \widetilde{M}$. Note also that for all elements $i$ of $[n]^\mc{Q}$, which satisfy (\ref{19.1.05.17.35}), the expectation  (\ref{21.1.05.1}) is the same, by invariance of Haar measure on unitary groups under permutation of rows and columns.
 Let us denote  this common expectation by $\E_\mc{Q}$. 
So,  one can write the sum of (\ref{emouvantgrosamericain'})\begin{equation}\label{emouvantgrosamericain'2005}
\ds
\E_\mc{Q}\sum_{i}\prod_{r=1}^{2N}D(r,n)_{ j_{M(r)}  }
\end{equation}
where the sum is taken on elements $i$ of $[n]^\mc{Q}$  which satisfy (\ref{19.1.05.17.35}). %\st for $h\in [2M]$ \st $i_h$ and $j_h$ are in the interval of integers $$\{q_1(n)+\cdots +q_{k(h)-1}(n)+1,\ldots,  q_1(n)+\cdots +q_{k(h)-1}(n)+q_{k(h)}(n)\}.$$
 
It is well known (see \cite{hiai}, equation (4.2.11)) that if, for each $m$, $U(m)$ is a uniform unitary $m\tii m$ random matrix, and $1\leq i_m,j_m \leq m$, then the sequence $\E (|U(m)_{i_m,j_m}|^{2M})$ does not depend on the choices of $i_m,j_m$ and \begin{equation}\label{pouleoeuf}\E (|U(m)_{i_m,j_m}|^{2M})^{1/\! 2M}=O(m^{-1/\! 2}).\end{equation}
So, by H\"older inequality, for 
$i\in [n]^\mc{Q}$, $$\lf|\E_\mc{Q}\ri|\leq\ds \prod_{h=1}^{2M}\E(|v_{i_h,j_h}(\kappa(h),\eps(h),n)|^{2M})^{1/\! 2M}=O(q_1(n)^{-M'}{n}^{-M''}),$$
where $\ds M'=\sum_{\substack{1\leq r\leq N\\ k_r=1}}|m_r|$ and $\ds M''=\sum_{\substack{1\leq r\leq N\\ k_r\neq 1}}|m_r|.$
Thus it suffices to prove that \begin{equation}\label{emouvantgrosamericain2005encoreplusgros}\ds\sum_{i}\prod_{r=1}^{2N}D(r,n)_{ j_{M(r)}  },\end{equation} is $o(q_k(n)^2q_1(n)^{M'}{n}^{M''})$ in case \uov and $o(q_1(n)^{M'+2}n^{M''})$ in case $2^\circ$),  where the sum is taken on the elements $i$ of $[n]^\mc{Q}$ which satisfy (\ref{19.1.05.17.35}).

{\it Step V. } \udeuzio We are going to prove it as an application of lemma \ref{01.11.03.1}.   The sum of (\ref{emouvantgrosamericain2005encoreplusgros}) is the sum, on the same $i$'s, of 
\begin{eqnarray*}&D(1,n)_{i_{M(1)+1}}D(2,n)_{i_{M(2)+1}}\cdots D(N-1,n)_{i_{M(N-1)+1}}D(N,n)_{i_1}&\\ \tii &D(N+1,n)_{i_{M+1}}D(N+2,n)_{i_{M(N+2)}}D(N+3,n)_{i_{M(N+3)}}\cdots D(2N,n)_{i_{M(2N)}}.&\end{eqnarray*}
i.e. of 
\begin{eqnarray*}&D(N,n)_{i_1}D(1,n)_{i_{M(1)+1}}D(2,n)_{i_{M(2)+1}}\cdots D(N-1,n)_{i_{M(N-1)+1}}&\\ \tii &D(N+1,n)_{i_{M+1}}D(N+2,n)_{i_{M(N+2)}}D(N+3,n)_{i_{M(N+3)}}\cdots D(2N,n)_{i_{M(2N)}}.&\end{eqnarray*} 
Adding some matrices of the type $p_k(n)$ in the product, one sees that this sum is equal to  \begin{equation}\label{emouvantgrosamericain2005encoreplusgrosetdansun4*4}  \ds\sum_{i\in [n]^\mc{Q}}\prod_{h=1}^{2M} A_{i_h}(h,n),\end{equation}
where for $h\in[2M]$, $$A(h,n)=\begin{cases}%p_{k(M)}(n)D(N,n)D(1,n)p_{k(1)}(n) &\textrm{if $h=1$}\\ 
p_{k(h-1)}(n)D(r,n)p_{k(h)}(n)&\textrm{if $\exists r\in\{0,\ldots,N-1\}$,}\\ &\quad h=M(r)+1\\
p_{k(h)}(n)&\textrm{if $h\in [M]$ and } h\notin\\ &\quad \{M(0)+1,\ldots,M(N-1)+1\}\\
p_{k(2M)}(n)D(N+1,n)p_{k(M+1)}(n)&\textrm{if $h=M+1$}\\
p_{k(h)}(n)D(r,n)p_{k(h+1)}(n)&\textrm{if $\exists r\in\{N+2,\ldots,2N\}$,}\\ 
 &\quad h=M(r)\\
p_{k(h)}(n)&\textrm{if $h\in [2N]$ and }h\notin \\ &\quad [M+1]\cup\{M(N+2),\ldots,M(2N)\}
\end{cases}$$
%Let       $\mc{H}$ be the set of elements $h$ \st $k(h)=1$. By definition of the $A(h,n)$'s, for all $h\in \mc{H}$, for all $i\in \{q_1(n)+1,\ldots,n\}$, $A_i(h,n)=0$.   %A DEPLACER !!! of $[M-1]$ \st $k(h)$ or $k(h+1)$ is $1$, plus $M$ if $k(M)$ or $k(1)$ is $1$, plus the set of elements $h$ of $\{M+2,\ldots, 2M\}$ \st $k(h-1)$ or $k(h)$ is $1$, plus $M+1$ if $k(M+1)$ or $k(2M)$ is $1$.

Now, note that  by lemma 4.2.2 of \cite{hiai} (or because the uniform distribution on the unitary group is left and right invariant), if $i\in [n]^\mc{Q}$  corresponds to a nonzero term in  equation (\ref{emouvantgrosamericain'}), then
for all $\kappa\in \n\tii [d]$,  for all $1\leq \alpha,\beta\leq n$, 
\begin{eqnarray*}
|\{h\ste i_h=\alpha,\kappa (h)=\kappa,\varepsilon (h)=1\}|&=&|\{h\ste j_h=\alpha,\kappa(h)=\kappa,\varepsilon (h)=-1\}|,\\
|\{h\ste j_h=\beta,\kappa (h)=\kappa,\varepsilon (h)=1\}|&=&|\{h\ste i_h=\beta,\kappa(h)=\kappa,\varepsilon (h)=-1\}|.
\end{eqnarray*}
Thus, if $\{h\}$ is a singleton class in $\mc{Q}$, and if $\E_\mc{Q}$ is non null, one has, using the cyclic orders on $[M]$ and on $\{M+1,M+2,\ldots, 2M\}$, 
\begin{eqnarray*}\varepsilon (h)=1 &\Rightarrow  & \kappa(h)=\kappa(h-1), \varepsilon (h-1)=-1,\\
{\textrm{and }}\qquad\varepsilon (h)=-1 &\Rightarrow  & \kappa(h)=\kappa(h-1), \varepsilon (h-1)=1.\end{eqnarray*}
It clearly follows that there exists $r\in\{0,\ldots, 2N\}-\{N+1\}$ \st $h=M(r)+1$. Moreover, by hypothesis, in case \uov   $A(h,n)=D(r,n)$ is a diagonal  matrix with null trace, and, in case $2^\circ$),  $A(h,n)=D(r,n)$ is a diagonal  matrix with null trace whenever $h\notin \{M,M+1\}$. 

\uo Let us apply the lemma with $\mc{P}=\mc{Q}$ and  $\mc{H}=\{h\in [2M]\ste k(h)=1\}.$ %being the set of elements $h$ \st $k(h)=1$. 
By definition of the $A(h,n)$'s, for all $h\in \mc{H}$, for all $i\in \{q_1(n)+1,\ldots,n\}$, $A_i(h,n)=0$. 
By the lemma, the sum of   (\ref{emouvantgrosamericain2005encoreplusgros}), which is equal to the sum of (\ref{emouvantgrosamericain2005encoreplusgrosetdansun4*4}), is $O(q_1(n)^{|\mc{P}_\mc{H}|-p_\mc{H}/2}n^{|\overline{\mc{P}}_\mc{H}|-\overline{p}_\mc{H}/2}).$
But for any partition $\mc{X}$ with $x$ singletons of a set $S$, one has $|\mc{X}|\leq x+(|S|-x)/2$, so $|\mc{X}|-x/2\leq |S|/2$. So, 
$$\ds |\mc{P}_\mc{H}|+|\overline{\mc{P}}_\mc{H}|-\f{p_\mc{H}}{2}-\f{\overline{p}_\mc{H}}{2}\leq (2M)/\!2 = M'+M'',$$
Moreover, since for all $h\in [2M]$,  $k(h)=1$ implies $h\in \mc{H}$, 
$\ds |\overline{\mc{P}}_\mc{H}|-\overline{p}_\mc{H}/2\leq 
\ff{2}\lf|\lf\{h\ste k(h)\neq 1\ri\}\ri|=M''$. %pour ça, il suffit d'avoir k(h)=1 => h in \mc{H}
Thus, since $q_1(n)\leq n$,  $n^{|\overline{\mc{P}}_\mc{H}|-\overline{p}_\mc{H}/2}=n^{|\overline{\mc{P}}_\mc{H}|-\overline{p}_\mc{H}/2-M''}n^{M''}=O(q_1(n)^{|\overline{\mc{P}}_\mc{H}|-\overline{p}_\mc{H}/2-M''}n^{M''})$. Hence$$ \ds q_1(n)^{|\mc{P}_\mc{H}|-p_\mc{H}/2}{n}^{|\overline{\mc{P}}_\mc{H}|-\overline{p}_\mc{H}/2} =O(q_1(n)^{|\mc{P}_\mc{H}|-\f{p_\mc{H}}{2}+|\overline{\mc{P}}_\mc{H}|-\f{\overline{p}_\mc{H}}{2}-M''}{n}^{M''})=O(q_1(n)^{M'}{n}^{M''}).$$
The sum of   (\ref{emouvantgrosamericain2005encoreplusgros}) is $O(q_1(n)^{M'}{n}^{M''})$, and the proposition is proved.

\deuzio Let us prove that the sum of (\ref{emouvantgrosamericain2005encoreplusgros}) is $o(q_1(n)^{2+M'}{n}^{M''})$ by induction on the number $|\mc{Q}|$ of classes of $\mc{Q}$. Let us denote, for $B, B'$ class of $\mc{Q}$, $\mc{Q}(M\to B)$ the partition of $[2M]$ obtained from $\mc{Q}$ by linking classes $B$ and $\{M\}$, and $\mc{Q}(M\to B, M+1\to B')$  the partition obtained from $\mc{Q}$ by linking classes $B$ and $\{M\}$, and then the classes $B'$ (which has become $B\cup\{M\}$ if $B'=B$) and $\{M+1\}$ Let $J(n)$ be the set of elements $j$ of $[n]$ \st $D(N,n)_j\neq 0$. Note that, as noted to step II, $|J(n)|\leq q_1(n),$ and, by definition of  $D(N+1,n)$, $J(n)$ is also the set of elements $j$ of $[n]$ \st $D(N+1,n)_j\neq 0$. 

$\bullet$ If $\mc{Q}$ has only one class, by boundedness hypothesis and since $|J(n)|\leq q_1(n)$, the sum of  (\ref{emouvantgrosamericain2005encoreplusgros}) is $O(q_1(n))$.% and the conclusion holds because, since $N\geq 1$, one have $M'\geq 1$ or $M''\geq 1$. 

$\bullet$ Assume that the conclusion holds to ranks $1,\ldots, |\mc{Q}|-1$.

- If  neither $\{M\}$ nor $\{M+1\}$ are singletons of $\mc{Q}$, we apply the lemma, as for case $1^\circ)$,  with $\mc{P}=\mc{Q}$ and   $\mc{H}=\{h\in [2M]\ste k(h)=1\}$.
It leads to the conclusion in the same way as in case \uop

- If exactly one of   $\{M\},\{M+1\}$, say $\{M\}$, is a singleton of $\mc{Q}$,   
define $\mc{P}=\mc{Q}-\{\{M\}\}$, partition of the set $[2M]-\{M\}$. Note that the sum of   (\ref{emouvantgrosamericain2005encoreplusgros}), which is equal to the sum of (\ref{emouvantgrosamericain2005encoreplusgrosetdansun4*4}), is  
\begin{eqnarray*} & \ds
\sum_{j\in J(n)} \Big[
 A(M,n)_j
\sum_{i\in ([n]-\{j\})^\mc{P}}\prod_{\substack{h=1\\ h\neq M}}^{2N} A_{i_h}(h,n)\Big]&\\  
=& \ds\sum_{j\in J(n)}\Big[A(M,n)_j\Big(\sum_{i\in [n]^\mc{P}}\prod_{\substack{h=1\\ h\neq M}}^{2N} A_{i_h}(h,n)-\sum_{\substack{i\in [n]^\mc{P}\\ \exists h, i_h=j}}\prod_{\substack{h=1\\ h\neq M}}^{2N} A_{i_h}(h,n)\Big)\Big]&\\
=&\ds\Big(\sum_{j\in J(n)}A(M,n)_j\Big)\Big(\sum_{i\in [n]^\mc{P}}\prod_{\substack{h=1\\ h\neq M}}^{2N} A_{i_h}(h,n)\Big)-\Big(\sum_{\substack{B\in \mc{Q}\\ B\neq \{M\}}}\sum_{i\in [n]^{\mc{Q}(M\to B}}\prod_{h=1}^{2N} A_{i_h}(h,n)\Big)&\\
\end{eqnarray*}
So by the induction hypothesis, for all $B\in \mc{P}, B\neq \{M\}$, we have $$\sum_{i\in [n]^{\mc{Q}(M\to B)}}\prod_{h=1}^{2N} A_{i_h}(h,n)=o\lf(q_1(n)^{2+M'}{n}^{M''}\ri).$$ So it suffices to prove that we have \begin{equation}\label{21.1.05.20.11}\Big(\sum_{j\in J(n)}A(M,n)_j\Big)\Big(\sum_{i\in [n]^\mc{P}}\prod_{\substack{h=1\\ h\neq M}}^{2N} A_{i_h}(h,n)\Big)=o\lf(q_1(n)^{2+M'}{n}^{M''}\ri).\end{equation} But by boundedness hypothesis, $\sum_{j\in J(n)}A(M,n)_j=0(|J(n)|)=0(q_1(n))$. So it suffices to prove that the second sum in (\ref{21.1.05.20.11}) is $o(q_1(n)^{M'+1}n^{M''})$.  
Define $\mc{H}=\{h\in [2M]-\{M\}\ste k(h)=1\}.$ 
By definition of the $A(h,n)$'s, for all $h\in \mc{H}$, for all $i\in \{q_1(n)+1,\ldots,n\}$, $A_i(h,n)=0$. 
So by the lemma,  the second sum in (\ref{21.1.05.20.11}) is $O(q_1(n)^{|\mc{P}_\mc{H}|-p_\mc{H}/2}n^{|\overline{\mc{P}}_\mc{H}|-\overline{p}_\mc{H}/2}).$  
Since, as noted above, for any partition $\mc{X}$ with $x$ singletons of a set $S$, one has $|\mc{X}|-x/2\leq |S|/2$, we have 
$$\ds |\mc{P}_\mc{H}|+|\overline{\mc{P}}_\mc{H}|-\f{p_\mc{H}}{2}-\f{\overline{p}_\mc{H}}{2}\leq (2M-1)/\!2 = M'+M''-1/\! 2,$$
Moreover, since for all $h\in [2M]-\{M\}$,  $k(h)=1$ implies $h\in \mc{H}$, 
$$\ds |\overline{\mc{P}}_\mc{H}|-\overline{p}_\mc{H}/2\leq 
\ff{2}\lf|\lf\{h\in [2M]-\{M\}\ste k(h)\neq 1\ri\}\ri|\leq M''.$$ %pour ça, il suffit d'avoir k(h)=1 => h in \mc{H}
Thus, since $q_1(n)\leq n$,  $n^{|\overline{\mc{P}}_\mc{H}|-\overline{p}_\mc{H}/2}=n^{|\overline{\mc{P}}_\mc{H}|-\overline{p}_\mc{H}/2-M''}n^{M''}=O(q_1(n)^{|\overline{\mc{P}}_\mc{H}|-\overline{p}_\mc{H}/2-M''}n^{M''})$. Hence $ \ds q_1(n)^{|\mc{P}_\mc{H}|-p_\mc{H}/2}{n}^{|\overline{\mc{P}}_\mc{H}|-\overline{p}_\mc{H}/2} =O(q_1(n)^{|\mc{P}_\mc{H}|-\f{p_\mc{H}}{2}+|\overline{\mc{P}}_\mc{H}|-\f{\overline{p}_\mc{H}}{2}-M''}{n}^{M''})=O(q_1(n)^{M'-1/\!2}{n}^{M''})$, thus the conclusion holds.
%The the sum of   (\ref{emouvantgrosamericain2005encoreplusgros}) is $O\lf(q_1(n)^{M'}{n}^{M''}\ri)$, and the proposition proved.

- If $\{M\}$ and $\{M+1\}$ are singletons  of $\mc{Q}$, define  $\mc{P}=\mc{Q}-\{\{M\}, \{M+1\}\}$, partition of the set $[2M]-\{M,M+1\}$.   Note that the sum of   (\ref{emouvantgrosamericain2005encoreplusgros}), which is equal to the sum of (\ref{emouvantgrosamericain2005encoreplusgrosetdansun4*4}), is  
%\begin{eqnarray*} & 
$$\ds
\sum_{\substack{j,j'\in J(n)\\ j\neq j'}} \Big[
 A(M,n)_jA(M+1,n)_{j'}
\sum_{i\in ([n]-\{j,j'\})^\mc{P}}\prod_{\substack{h=1\\ h\neq M\\ h\neq M+1}}^{2N} A_{i_h}(h,n)\Big].$$

Now, let us define $\mc{Q}'$ to be the partition of $[2M]$ obtained from $\mc{Q}$ by linking classes $\{M\}$ and $\{M+1\}$.
The sum of (\ref{emouvantgrosamericain2005encoreplusgros}) can be written$$\ds
\Big(\sum_{j,j'\in J(n)} 
 A(M,n)_jA(M+1,n)_{j'}\Big)\Big(
\sum_{i\in [n]^\mc{P}}\prod_{\substack{h=1\\ h\neq M\\ h\neq M+1}}^{2N} A_{i_h}(h,n)\Big)-\sum_{i\in [n]^{\mc{Q}'}}\prod_{h=1}^{2N} A_{i_h}(h,n)-\sum_{\mc{R}}\sum_{i\in [n]^{\mc{R}}}\prod_{h=1}^{2N} A_{i_h}(h,n),$$ where in the last sum, $\mc{R}$ runs over  the set $\{\mc{Q}(M\to B,M+1\to B')\ste B,B'\in \mc{Q}\}$.
So by the induction hypothesis, it suffices to prove that  \begin{equation}\label{23.1.05.1}\Big(\sum_{j,j'\in J(n)} 
 A(M,n)_jA(M+1,n)_{j'}\Big)\Big(
\sum_{i\in [n]^\mc{P}}\prod_{\substack{h=1\\ h\neq M\\ h\neq M+1}}^{2N} A_{i_h}(h,n)\Big)=o\lf(q_1(n)^{2+M'}{n}^{M''}\ri).\end{equation}
But by boundedness hypothesis, $\sum_{j,j'\in J(n)}A(M,n)_jA(M+1,n)_{j'}=0(|J(n)|^2)=0(q_1(n)^2)$. So it suffices to prove that the second sum in (\ref{23.1.05.1}) is $o(q_1(n)^{M'}n^{M''})$.  
Define $\mc{H}:=\{h\in [2M]-\{M,M+1\}\ste k(h)=1\}$. By definition of the $A(h,n)$'s, for all $h\in \mc{H}$, for all $i\in \{q_1(n)+1,\ldots,n\}$, $A_i(h,n)=0$. So by the lemma,  the second sum in (\ref{23.1.05.1}) is $$O(q_1(n)^{|\mc{P}_\mc{H}|-p_\mc{H}/2}n^{|\overline{\mc{P}}_\mc{H}|-\overline{p}_\mc{H}/2}).$$ Since, as noted above, for any partition $\mc{X}$ with $x$ singletons of a set $S$, one has $|\mc{X}|-x/2\leq |S|/2$, we have 
$$\ds |\mc{P}_\mc{H}|+|\overline{\mc{P}}_\mc{H}|-\f{p_\mc{H}}{2}-\f{\overline{p}_\mc{H}}{2}\leq (2M-2)/\!2 = M'+M''-1,$$ Moreover, since for all $h\in [2M]-\{M,M+1\}$,  $k(h)=1$ implies $h\in \mc{H}$, 
$$\ds |\overline{\mc{P}}_\mc{H}|-\overline{p}_\mc{H}/2\leq 
\ff{2}\lf|\lf\{h\in [2M]-\{M,M+1\}\ste k(h)\neq 1\ri\}\ri|\leq M''.$$ %pour ça, il suffit d'avoir k(h)=1 => h in \mc{H}
Thus, since $q_1(n)\leq n$,  $$n^{|\overline{\mc{P}}_\mc{H}|-\overline{p}_\mc{H}/2}=n^{|\overline{\mc{P}}_\mc{H}|-\overline{p}_\mc{H}/2-M''}n^{M''}=O(q_1(n)^{|\overline{\mc{P}}_\mc{H}|-\overline{p}_\mc{H}/2-M''}n^{M''}).$$ Hence$$ \ds q_1(n)^{|\mc{P}_\mc{H}|-p_\mc{H}/2}{n}^{|\overline{\mc{P}}_\mc{H}|-\overline{p}_\mc{H}/2} =O(q_1(n)^{|\mc{P}_\mc{H}|-\f{p_\mc{H}}{2}+|\overline{\mc{P}}_\mc{H}|-\f{\overline{p}_\mc{H}}{2}-M''}{n}^{M''})=O(q_1(n)^{M'-1}{n}^{M''}),$$ thus the conclusion holds.
\end{pr}

 \end{document}